%%%%%%%%%%%%%%%%%%%%%%%%%%%%%%%%%%%%%%%%%%%%%%%%    
%
%        THIS IS A  PLAIN TeX FILE
%
%%%%%%%%%%%%%%%%%%%%%%%%%%%%%%%%%%%%%%%%%%%%%%%%

\magnification=1200

\font\titfont=cmr10 at 12 pt

\font\headfont=cmr10 at 12 pt

%\font\AAA=Times.dfont  at 12pt
 %\font\BBB=Times.dfont at 8pt

%\font\AAA=cmr10 at 12pt
%\font\BBB=cmr10 at 8pt

\def\AAA{\bf}
\def\BBB{\bf}

\overfullrule=0in

\def\boxit#1{\hbox{\vrule
 \vtop{%
  \vbox{\hrule\kern 2pt %
     \hbox{\kern 2pt #1\kern 2pt}}%
   \kern 2pt \hrule }%
  \vrule}}
  
\def\imp{\qquad \Rightarrow\qquad}

  \def\harr#1#2{\ \smash{\mathop{\hbox to .3in{\rightarrowfill}}\limits^{\scriptstyle#1}_{\scriptstyle#2}}\ }

\def\intave#1{\int_{#1}\!\!\!\!\! \! \!\!\!\!-}

 \def\GG{{{\bf G} \!\!\!\! {\rm l}}\ }

\def\snm{{S^{n-1}}}

\def\ST{{ST-invariant}}

\def\snm{S^{n-1}}
\def\den{\Theta}
\def\sh{subharmonic}
\def\iiff{\qquad\iff\qquad}
\def\bra#1#2{\langle #1, #2\rangle}

\def\ss{\subset}

\def\smfrac#1#2{\hbox{${#1\over #2}$}}

\def\dim{{\rm dim}}

\def\log{{\rm log}}
\def\Hess{{\rm Hess}}

\def\tr{{\rm tr}}

\def\span{{\rm span\,}}

\def\det{{\rm det}}

\def\Sym{{\rm Sym}^2}

\def\arr{\longrightarrow}

\def\rn{\bbr^n}

\def\Int{{\rm Int}}

\def\Symn{{\Sym(\rn)}}

\def\Theorem#1{\medskip\noindent {\bf THEOREM \bf #1.}}
\def\Prop#1{\medskip\noindent {\bf Proposition #1.}}
\def\Cor#1{\medskip\noindent {\bf Corollary #1.}}
\def\Lemma#1{\medskip\noindent {\bf Lemma #1.}}
\def\Remark#1{\medskip\noindent {\bf Remark #1.}}
\def\Note#1{\medskip\noindent {\bf Note #1.}}
\def\Def#1{\medskip\noindent {\bf Definition #1.}}

\def\Ex#1{\medskip\noindent {\bf Example \bf    #1.}}

\def\pf{\medskip\noindent {\bf Proof.}\ }
\def\qed{\hfill  $\vrule width5pt height5pt depth0pt$}

\def\cv{{\cal V}}   \def\cp{{\cal P}}

\def\ch{{\cal H}}

\def\cl{{\cal L}}
\def\cp{{\cal P}}

\def\vf{\varphi}

\def\wt{\widetilde}

\def\and{\qquad {\rm and} \qquad}
\def\arr{\longrightarrow}
\def\ol{\overline}
\def\bbr{{\bf R}}\def\bbh{{\bf H}}
\def\bbc{{\bf C}}

\def\bbp{{\bf P}}

\def\bbe{{\bf E}}

\def\d{\delta}

\def\g{\gamma}

\def\l{\lambda}
\def\o{\omega}

\def\s{\sigma}

\def\D{\Delta}
\def\L{\Lambda}
\def\G{\Gamma}
\def\O{\Omega}

\def\lloc{L^1_{\rm loc}}
\def\dbar{\ol{\partial}}

\def\Symn{\Sym(\rn)}
 
\def\USC{{\rm USC}}
\def\fa{{\rm\ \  for\ all\ }}

\def\cpt{\wt{\cp}}
\def\ft{\wt F}

 \def\LAG{{\rm LAG}}
 \def\ISO{{\rm ISO_p}}

\def\AA{1}
\def\BB{2}
\def\CC{3}

\def\EE{5}
\def\FF{6}

\def\HH{7}

\def\JJ{10}
 \def\KK{11}
 \def\LL{12}
\def\MM{13}

\def\HLCG{HL$_1$} 
\def\HLPCG{HL$_2$} 
 
\def\HLDD{HL$_3$} 
 
\def\HLDDR{HL$_4$}
 
\def\HLHPP{HL$_5$} 
\def\HLREST{HL$_6$} 
\def\HLPCON{HL$_{7}$} 
\def\HLSURVEY{HL$_{8}$}

\def\HLRS{HL$_{9}$}
\def\HLTANI{HL$_{10}$}

\def\intave#1{\int_{#1}\!\!\!\!\!\!\!-\ }

\def\AAA{1}
 \def\AA{2}
 \def\CC{3}
 \def\HH{4} 
 \def\BB{5}
 \def\BBB{6}
  \def\FF{7} 
  \def\FFF{8}
\def\JJ{9} 
\def\JJJ{10}
 \def\MM{11}
  
\def\GGG{2}
\def\KK{3}  
\def\EE{4} 
\def\LL{5}

\def\ST{ST-invariant }

\def\point{x}

\centerline
{
\titfont  TANGENTS TO SUBSOLUTIONS}
\medskip

\centerline{\titfont EXISTENCE AND UNIQUENESS, II}

\bigskip

\centerline{\titfont F. Reese Harvey and H. Blaine Lawson, Jr.$^*$}
\vglue .9cm
\smallbreak\footnote{}{ $ {} \sp{ *}{\rm Partially}$  supported by
the N.S.F. }

%\vskip .1in
 
\centerline{\bf ABSTRACT} \medskip
 % \font\abstractfont=cmr10 at 10 pt

  {{\parindent= .44in
\narrower 
This part II of the paper is concerned with questions of existence
and uniqueness of tangents in the special case of  $\GG$-plurisubharmonic functions,
where $\GG \ss G(p,\rn)$ is a  compact subset of the 
Grassmannian of $p$-planes in $\rn$.  An u.s.c. function $u$ on an
open set $\O\ss\rn$ is $\GG$-plurisubharmonic if its restriction to
$\O\cap W$ is subharmonic for every affine $\GG$-plane $W$.
Here $\GG$ is assumed to be invariant under a subgroup $K\ss {\rm O}(n)$
which acts transitively on $S^{n-1}$.
Tangents to $u$ at a point $\point$ are the cluster points of $u$ under 
a natural flow (or blow-up) at $\point$. 
They always exist and are $\GG$-harmonic at all points of continuity.
A homogeneity property is established for all tangents in these geometric cases.
This leads to principal results concerning the Strong Uniqueness of Tangents, which
means that all tangents are unique and of the form $\den K_p$ where
$K_p$ is the Riesz kernel and $\den$ is the density of $u$ at the point.
Strong uniqueness is a form of regularity
which implies that the sets $\{\den(u,x)\geq c\}$ for $c>0$ are discrete.
When the invariance group $K= {\rm O}(n), {\rm U}(n)$ or Sp$(n)$ strong uniqueness 
holds for all but a small handful of cases.
It also holds for essentially  all interesting $\GG$ which arise in calibrated geometry.

When strong uniqueness fails, homogeneity implies that tangents are characterized
by a subequation on the sphere, which is worked out in detail.
In the cases corresponding to the real, complex and quaternionic Monge-Amp\`ere
equations (convex functions, and complex and quaternionic plurisubharmonic
functions) tangents, which are far from unique, are then systematically studied and
classified.

}}

\vskip.3in

\centerline{\bf TABLE OF CONTENTS} \bigskip
% \font\abstractfont=cmr10 at 10 pt

{{\parindent= .1in\narrower  \noindent

\qquad \AAA. Introduction.\smallskip

\qquad\GGG.  The Homogeneity Theorem.\smallskip

\qquad\KK.  The Strong Uniqueness  Theorem.\smallskip

\qquad \EE.    Homogeneous F-Subharmonic Functions.     \smallskip

\qquad \LL. Tangents to Convex, $\bbc$-Plurisubharmonic,  and $\bbh$-Plurisubharmonic

\qquad \quad\  \ \ Functions.

\vskip.15in

\hskip.5 in
Appendix A.    Further Discussion of Examples.

}}

\vfill\eject

%%%%%%%%%%%%%%%%%%%%%%%%%%%%%%%%%%%%%%%%%%%%%%%%%%%%
%%%%%%%%%%%%%%%%%%%%%%%%%%%%%%%%%%%%%%%%%%%%%%%%%%%%
%%%%%%%%%%%%%%%%%%%%%%%%%%%%%%%%%%%%%%%%%%%%%%%%%%%%
%%%%%%%%%%%%%%%%%%%%%%%%%%%%%%%%%%%%%%%%%%%%%%%%%%%%
%%%%%%%%%%%%%%%%%%%%%%%%%%%%%%%%%%%%%%%%%%%%%%%%%%%%

\centerline{\headfont \AAA.   Introduction.}
\medskip 
Part I of this paper was concerned with the study of tangents to $F$-subharmonic functions
(or  subsolutions) for any fully nonlinear subequation $F\ss \Symn$. 
Key to the results is the notion of the {\sl Riesz characteristic of $F$}, a real number
 $p=p_F$ with $1\leq p\leq \infty$.   When  $p$  is finite,  there is  an associated tangential $p$-flow
on $F$-subharmonic functions $u$  at %each point $x_0$ given (for $x_0=0$) by
0 given by
$$
u_r(x) \ = \ 
\cases
{
\ \ r^{p-2}  u(rx) \qquad\qquad {\rm if}\ \ p\neq 2, \  \ {\rm and}   \cr
  u(rx)  -  M(u,r)    \qquad {\rm if}\ \ p=2, 
}
\eqno{(\AAA.1)}
$$
where 
$$
M(r) \ \equiv \  \sup_{|x|\leq r} u.
\eqno{(\AAA.2)}
$$
 Tangents to $u$ at 0 are then defined to be the 
cluster points of this flow in $\lloc(\rn)$. A basic result is that {\bf tangents always exist},
and the set of tangents to $u$ at 0 has a list of  characterizing properties 
(Part I, Section \FFF). Tangents are also always maximal (Part I, Section \BBB).
In particular, they are $F$-harmonic outside possible poles.

Of basic importance to this study is the {\bf  $p^{\rm th}$ Riesz kernel} $K_p(|x|)$ where
$$
K_p(t)\ =\ \cases
{
t^{2-p}          \qquad {\sl if} \ \ 1\leq p\ <2             \cr
\log \, t          \qquad {\sl if} \ \ p=2        \cr
-{1\over t^{p-2}}          \quad {\sl if} \ \  2<p<\infty.           \cr
}
\eqno{(\AAA.3)}
$$
When the Riesz characteristic $p=p_F$ is finite, every increasing radial $F$-harmonic
is of the form $\den K_p(|x|) +C$.  A fundamental Monotonicity Theorem (Part I, Section \BB)
states that 
$$
{M(u,r)-M(u,s)  \over K_p(r)-K_p(s)} \quad {\rm is\ increasing\ in} \ \ r\ \ {\rm and}\ \ s.
\eqno{(\AAA.4)}
$$
for all $0<r<s$ where $M$ is defined.  This gives the notion of the {\bf density
of $u$ at $0$}:
$$
\den(u,0) \ =\ \lim_{r<s\downarrow0} {M(u,r)-M(u,s)  \over K_p(r)-K_p(s)}.
\eqno{(\AAA.5)}
$$
(When $F$ is convex, there are other densities defined via the  area and 
volume averages.)

This part of the paper is exclusively concerned with {\sl geometric subequations}
determined by a closed subset $\GG\ss G(p,\rn)$ of the Grassmannian of 
unoriented $p$-planes in $\rn$  (where $1\leq p\leq n$).
 We recall that the associated subequation is
$$
F(\GG) \ \equiv\ \left \{A\in \Symn : \tr\left(A\bigr|_W\right)\geq0 \ \forall\, W\in \GG\right\}.
\eqno{(\AAA.6)}
$$
The Riesz characteristic of   $F(\GG)$ is the integer $p$.
The $F(\GG)$-subharmonic functions are called  {\sl $\GG$-plurisubharmonic}, 
and they are characterized by the property that their restrictions to   affine $\GG$-planes 
are subharmonic [\HLREST]. Many examples of geometric interest are given in Part I.  These include
in particular the plurisubharmonics associated to any calibration, and the Lagrangian
subharmonics in $\bbc^n$.

Recall that the standing assumptions on  $F\ss \Symn$ in Part I were:\ 

\qquad (i) \ \ (Positivity) $F+\cp\ss F$ where $\cp=\{A\geq0\}$, \ 

\qquad  (ii)\ \  (Cone property) $tF=F$ for all $t\geq0$,

\qquad  (iii)\  (ST-Invariance) $F$ is invariant under a subgroup $G\ss {\rm O}(n)$ 

\qquad \qquad\qquad which acts transitively
on the sphere $S^{n-1}$. 

\noindent
The first two assumptions are automatic for $F(\GG)$.  The last is 
equivalent to the assumption that 
$\GG$ is invariant under the subgroup $G\ss {\rm O}(n)$ acting on the Grassmannian $G(p,\rn)$.

To state the main results we recall the following.

\Def{\AAA.1} We say that {\bf uniqueness of tangents} holds for the subequation $F$
if for every  $F$-subharmonic function $u$  defined in a neighborhood of $0$, 
there is exactly one tangent to $u$ at 0.
We say that {\bf strong uniqueness of tangents} holds for  $F$
if for every such $u$,  the unique tangent is $\den(u,0) K_p(|x|)$.

\Def{\AAA.2} An upper semi-continuous function $U:\rn\to [-\infty, \infty)$ is 
said to have {\bf Riesz homogeneity $p$} if $U_r = U$ for all $r>0$. 
This condition holds if and only if 
there exists an u.s.c. function $g$ on the
unit sphere $S$ such that 
$$
U(x) \ =\ |x|^{p-2} g\left(  {x\over |x|}\right) \qquad{\rm in\ the\ cases\  where \ \ } p\neq2,
\eqno{(\AAA.7)}
$$
while in the case where $p=2$, 
$$
U(x)\ =\ \den \log |x| + g\left(  {x\over |x|}\right) \qquad {\rm with}\ \ \sup_{S^{n-1}} g \ =\ 0\ \ \ {\rm and}\ \ \ 
\den\geq 0\ \ {\rm a \ constant.}
\eqno{(\AAA.8)}
$$

\medskip

\Note {\AAA.3}  When $p=1$ our assumption of ST-invariance implies
that $\GG = G(1,\rn)$.  Hence, there is only one geometric subequation, namely $\cp = \{A\geq0\}$.
The $\cp$-subharmonic functions are exactly the  convex functions, and
in this case straightforward classical arguments establish the existence, uniqueness and homogeneity of 
tangents at every point.  These proofs are omitted. On the other hand, strong uniqueness of tangents
fails in this case, and the classification is given in section \LL.
\medskip

Our first main result is the following.  Let $\GG \ss G(p,\rn)$ be as above, with $p\geq2$.

\medskip
\noindent
{\bf THE HOMOGENEITY THEOREM.} {\sl
Suppose $u$ is a $\GG$-plurisubharmonic function defined in a neighborhood of
0, and suppose $U$ is a tangent to $u$ at 0.
Then $U$ has Riesz homogeneity $p$.  
Moreover, for all $\GG$-planes $W$ passing through the origin, 
the function $g$ is constant on the unit sphere $W\cap S^{n-1}$ in $W$.
In fact, when $p\neq 2$, }
$$
g\left({x \over |x|}\right) \ =\ -\den(W) \qquad {\sl if}\ \  x\in W\in \GG
$$

The fact 
that  $g$   is constant on each intersection $W\cap S^{n-1}$
for $W\in\GG$ leads to the following.  We say $\GG$ has the  {\bf transitivity property} if 
 for any two  vectors $x,y\in\rn$ there exist $W_1,...,W_k\in\GG$
with $x\in W_1, y \in W_k$ and dim$(W_i\cap W_{i+1}) >0$ for all $i=1,...,k-1$.

\medskip
\noindent
{\bf  THE PRELIMINARY  STRONG UNIQUENESS THEOREM.} {\sl
If $\GG$ has the transitivity property, then strong uniqueness of tangents holds
for all $\GG$-plurisubharmonic functions.}

\medskip

This covers a number of interesting cases which are not included in the Strong Uniqueness Theorem 
of Part I.  For example, this establishes strong uniqueness for plurisubharmonic functions in
Special Lagrangian, associative, coassociative and Cayley geometry (See Section \KK.)
The invariance groups  in these cases are
 SU$(n)$, G$_2$ and Spin$_7$. For the standard families of groups acting transitively
 on spheres we have the following nearly complete result.

\medskip
\noindent
{\bf  THE PRINCIPAL STRONG UNIQUENESS THEOREM.} {\sl
  Fix $p\geq 2$ and $n\geq 3$.
Then strong uniqueness of tangents  to  $\GG$-plurisubharmonic functions
holds for:
\smallskip

(a) Every compact SU$(n)$-invariant subset 
$
\GG \ss G^\bbr(p, \bbc^n)
$
except $\cp^\bbc$,
\smallskip

(b) Every compact Sp$(n)\cdot {\rm Sp}(1)$-invariant subset 
$
\GG \ss G^\bbr(p, \bbh^n)
$
with three exceptions, namely 
 the sets of real $p$-planes which lie in 
 a quaternion line for $p=2,3,4$  (when $p=4$ this is $\cp^\bbh$),

(c) Every compact Sp$(n)$-invariant subset 
$
\GG \ss G^\bbr(p, \bbh^n),
$
 for $p\geq5$.
}

\medskip

We recall that by  Theorem \MM.1 from Part I: {\sl If strong uniqueness holds, 
then for every $\GG$-plurisubharmonic function $u$, the set }
$$
E_c \ =\ \left \{ x : \den(u,x)\geq c \right \} \ \ {\sl is\ discrete \ for\ all}\ \ c>0.
$$

In those cases where strong uniqueness fails  we have the following question :  
What is the subequation on the 
sphere $S^{n-1}$ satisfied by the function $g$ in (\AAA.7)?
This subequation is worked out in Section \EE.   Its viscosity subsolutions  are
exactly the functions $g$ in the Homogeneity Theorem above.

The three classical cases where strong uniqueness of tangents fails are: 
$$
\GG \ =\ G^\bbr(1, \rn), \ \  G^\bbc(1, \bbc^n) \ \ {\rm and}\ \ G^\bbh(1, \bbh^n).
$$
The associated subequations
$$
F(\GG) \ =\ \cp^\bbr, \ \  \cp^\bbc \ \ {\rm and}\ \ \cp^\bbh,
$$
the homogeneous real, complex and quaternionic Monge-Amp\`ere equations
respectively.  These cases are discussed in detail in Section \LL.

For the first case, $\cp^\bbr$-subharmonic functions are just classical convex 
functions.  Here tangents are unique, but strong uniqueness is far from true.
The results here are classical, but we review them for the light they shed on the
general picture.

For the second case, $\cp^\bbc$-subharmonic functions are the standard
plurisubharmonic functions in $\bbc^n$.  Here even the uniqueness of tangents fails.
However, the subsets of functions in $\lloc(\bbc^n)$ which can arise as the 
set of tangents at 0 to a p.s.h. function $u$ have been completely classified by
Kiselman [K] whose work was the inspiration for this paper. One new feature of
our presentation is that in this case, we show that  tangents correspond 
bijectively  to  quasi-plurisubharmonic
functions on complex projective space $\bbp^{n-1}(\bbc)$.

For the third case, $\cp^\bbh$-subharmonic functions are quaternionic plurisubharmonic 
functions (cf. [A$_{1,2}$], [AV]). Here the determination of tangents is new.
As above, the  tangents corresponds bijectively to upper semi-continuous 
functions $g$ on quaternionic projective space $\bbp^{n-1}(\bbh)$ which satisfy the inequality
$\Hess_\bbh (g) -2gI\geq0$ in the viscosity sense.

Finally, in  Appendix A we give a rounded discussion of the many examples to 
which the results of both Parts I and II apply.  This includes the  establishment of the maximal and minimal subequations of Riesz characteristic $p$ as well as the maximal and minimal ones 
in the convex case.

\vskip.3in
%\vfill\eject

%%%%%%%%%%%%%%%%%%%%%%%%%%%%%%%%%%%%%%%%%%%%%%%%%%%%
%%%%%%%%%%%%%%%%%%%%%%%%%%%%%%%%%%%%%%%%%%%%%%%%%%%%
%%%%%%%%%%%%%%%%%%%%%%%%%%%%%%%%%%%%%%%%%%%%%%%%%%%%
%%%%%%%%%%%%%%%%%%%%%%%%%%%%%%%%%%%%%%%%%%%%%%%%%%%%
%%%%%%%%%%%%%%%%%%%%%%%%%%%%%%%%%%%%%%%%%%%%%%%%%%%%

\centerline{\headfont \GGG.    The Homogeneity Theorem.}
\medskip 

In this section we establish the homogeneity of tangents for all  geometrically determined subequations.
We assume,  to begin, that  $\GG\ss G(p,\rn)$  is a smooth compact submanifold of the Grassmannian
where the integer $p$ equals $2, 3, ..., n-1$. Later we will be able to drop this smoothness assumption
and allow $\GG$ to be any closed subset of $p$-planes.
We always assume that $\GG$ is invariant under the natural action on  $G(p,\rn)$ of a subgroup $G\ss {\rm O}(n)$ which acts transitively on the unit sphere $S^{n-1}\ss \rn$.
 This implies in particular that every vector in $\rn$ lies in some $\GG$-plane

Recall that the associated subequation  $F$ is given by 
$$
F\ =\ F(\GG) \ \equiv\ \{A\in\Symn : \tr_W A\geq 0 \fa W\in \GG\}.
$$
The Riesz characteristic of $F$ is easily seen to be the integer $p$, 
and in fact, this is the reason for the choice of normalization for the Riesz kernel
$K_p$ in (\AAA.3).

Suppose $u$  is an $F$-\sh\  function which is defined in a neighborhood of the origin,
and $U$ is a tangent to $u$ at 0.  We assume $u\not\equiv -\infty$.
By the Restriction Theorem proved in [\HLREST]
$$
u\bigr|_W \ \ {\rm is\ Laplacian\ subharmonic\  on \ } W \ {\rm (near\ 0)\ for \ each \ } W
\in \GG.
\eqno{(\GGG.1)}
$$
In particular, either $u\bigr|_W$ is $\lloc$ or $u\bigr|_W \equiv -\infty$.
We say that $W$ is {\sl non-polar} for $u$ at 0 in the first case, and 
{\sl polar}   in the second case.
The invariance of $\GG$ implies that $F\ss \D$ (see (\BB.3) in Part I), and hence $u$ is $\D$-subharmonic
near $0\in\rn$.  Therefore, its $-\infty$ set has $\D$-capacity zero, and hence measure zero.
This proves
$$
{\rm The \ union\  of \ all\  polar\  planes\ } W\in\GG\ {\rm is \ a \ set\  of \ measure\  zero\ in\ } \rn
\eqno{(\GGG.2)}
$$
If $W$ is  non-polar for $u$ near 0, then because of (\GGG.1) we can apply the classical fact for the 
Laplacian that
$$
u\bigr|_W \ \ 
{\rm has\ the \ unique\ tangent\ function\ \  } \den  (W)K_p
\eqno{(\GGG.3)}
$$
where $K_p$ denote the function $K_p(|x|)$ and  where the constant 
$$
\den (W) = \den^M \left(u\bigr|_W\right)= \den^S \left(u\bigr|_W\right)
$$
is  the maximum and/or spherical density of $u\bigr|_W$ at 0  (see (\JJ.3) and Proposition \JJ.4 in 
Part I).  That is,
 $$
 \lim_{r\to 0} \left ( u \bigr|_W \right )_r   \ =\  \den\left (W\right) K_p
\quad\  \ \qquad {\rm in\ \ } \lloc(W) 
 \eqno{(\GGG.3)'}
$$
Note that  these limits are  over {\bf all} $r$, not just a sequence $r_j$.
Also note that    for $p\geq3$ we have 
$$
 \left ( u \bigr|_W \right )_r =  u_r \bigr|_W.
  \eqno{(\GGG.4)}
$$
%while for $p=2$,   $\left ( u \bigr|_W \right )_r  = u_r  \bigr|_W + M(u,r) - M( u \bigr|_W, r)$
(This does not hold for $p=2$.)

The main result of this section is the following. 
Recall we assume $u\not\equiv -\infty$.
Extend the definition of $\den(W)$ to all $W\in \GG$ by defining $\den(W) = +\infty$ 
if $W$ is polar for $u$ at 0.

\Theorem {\GGG.1}  {\sl
Suppose $u$  is an $F$-\sh\  function which is defined in a neighborhood of the origin,
and $U$ is a tangent to $u$ at 0.
Then  $U$ has Riesz homogeneity $p$,  that is,
$$
U(x) \ =\ 
\cases
{
 \  \ {1 \over |x|^{p-2}}\,  g\bigl({x \over |x|}\bigr) \qquad \ \ \  {\rm if}   \ \ p> 2  \cr \cr
\den \log\, |x| + \,g\bigl({x \over |x|}\bigr)  \qquad {\rm if} \ \ p= 2
}
\qquad {\sl where}\ \ g\ \equiv \ U\bigr|_{S^{n-1}}  \in  \USC(S^{n-1}).
 \eqno{(\GGG.5)}
$$ 
and where, in the case $p=2$, $\sup g =0$ and $\den=\den^M(u,0)$.  Moreover, 
for each $\GG$-plane $W$ passing through the origin, the function $g$ is constant 
on the unit sphere $W\cap S^{n-1}$ in $W$. In fact, when $p > 2$,
}
%$$
%{\sl if}\ \ x\in W \in \GG, \ \ {\sl then}\ \ g\left( {x\over |x|} \right)  \ =\ - \den(W).   
% \eqno{(\GGG.6)}
%$$ 
$$
g\left( {x\over |x|} \right)  \ =\ - \den(W) \qquad {\rm for} \ \ x\in W \in\GG.   
 \eqno{(\GGG.6)}
$$

\pf   \def\An{{\bf A}}
We first treat the case $p\geq3$.
Fix constants $0<a<b$, and let 
$\An = \{x\in \rn: a\leq |x|\leq b\}$ be the annulus with radii $a,b$.
 It will suffice to prove  our assertions  on  $\An$.

To begin set $Gr \equiv G(p,\rn)$ and  consider the tautological vector bundle
$$
E\  \equiv\ \{(W,x) \in Gr\times \rn : x\in W\}\  \harr{\s}\ \ Gr
$$
where $\s$ is given by projection onto the first factor in $Gr \times \rn$.
Projection onto the second factor gives another map
$$
\matrix
{
\ &\ & E & \ &\ \cr
\cr
\ \ \ \sigma &\swarrow & \ & \searrow & \pi  \cr
\cr
Gr&\ &\ &\ &  \rn
}
$$ 
Note that $\pi : E-Z \to \rn-\{0\}$ is a proper submersion, where 
$Z=\pi^{-1}(0) \ss Y$ is the zero section of the vector bundle $\s$.
Setting $E_\An = \pi^{-1}(\An)$ we have a pair of smooth compact fibre bundles
$$
\matrix
{
\ &\ & E_\An & \ &\ \cr
\cr
\sigma &\swarrow & \ & \searrow & \pi  \cr
\cr
Gr &\ &\ &\ &  \An
}
$$ 
where the fibre of $\s$ over $W\in Gr$ is the $(a,b)$-annulus in $W$. 
Note that the orthogonal group acts naturally on this diagram.

We now restrict this annulus-bundle $E_\An$ to the submanifold $\GG\ss Gr$, that is, we set
$\bbe \equiv \s^{-1}(\GG)$.  The diagram above reduces to a new diagram
$$
\matrix
{
\ &\ & \bbe & \ &\ \cr
\cr
\sigma &\swarrow & \ & \searrow & \pi  \cr
\cr
\GG &\ &\ &\ &  \An
}
\eqno{(\GGG.7)}
$$ 
Note that the subgroup $G\ss {\rm O}(n)$ acts naturally on this diagram (\GGG.7),
and recall that $G$ acts transitively on the concentric spheres $S_r^{n-1} = \{|x|=r\}, \  a\leq r\leq b$, in the
annular region $\An$. This, together with the fact that $\pi$ is a linear embedding on the fibres of 
$\s$, shows that $\bbe \harr{\pi}\  \An$ is also a smooth fibre-bundle 
over the manifold-with-boundary $\An$.

We fix a defining sequence
 $u^j\equiv u_{r_j}$ for $U$,  and  consider 
 the pull-backs $\wt{u}^j \equiv  \pi^* u^j$ and $\wt U \equiv \pi^* U$ to $\bbe$.
 Note that $\wt U$ is u.s.c., in fact it is essentially u.s.c. (since $U$ is), and 
 we have that
$$
\wt{u}^j  \ \to\ \wt U\qquad{\rm in}\ \ L^1(\bbe).
\eqno{(\GGG.8)}
$$

In addition, set $\wt V(W,x) \equiv \den(W)K_p(|x|)$ if $x\neq 0$ and $(W,x)\in \bbe$, i.e., $x\in W$.
Then (\GGG.3)$'$ implies that
$$
\wt{u}^j \bigr|_{\s^{-1}(W)}(x)  \ \to\   \wt V(W,x)
\quad{\rm in}\ \ L^1(\s^{-1}(W)))  \cong  L^1(W\cap \An)) 
\quad \forall\, W\in\GG\ \ {\rm non polar}.
\eqno{(\GGG.9)}
$$

\Lemma {\GGG.2} {\sl  $ \wt  U =  \wt  V$ almost everywhere on $\bbe$.
Furthermore,}
$$
U\bigr|_W \ =\  \den(W) K_p \qquad {\sl for \ almost \ all \ } W\in\GG.
\eqno{(\GGG.10)}
$$

\Cor{\GGG.3} {\sl $U_r=U$ on $\rn$ for all $r>0$, i.e., $U$ is $p$-homogeneous.}

\medskip
\noindent
{\bf Proof of Corollary \GGG.3}.
By (\GGG.10) and (\GGG.2) we see that   $U=U_r$ a.e. in $\rn$. 
Note, however,   that $U=U_r$ a.e. implies that $U=U_r$ everywhere
since both functions are classically $\D$-subharmonic
(and therefore satisfy $U(x) = \lim_{r\to 0} {\rm ess\, sup}_{B_r(x)} U$ for all $x$).
 \qed  \medskip

\medskip
\noindent
{\bf Proof of Lemma \GGG.2}.
The fibration $\s:\bbe\to \GG$ is locally a product $B\times A$ where $B$ is an open ball
in the manifold $\GG$ and $A$ is the $[a,b]$-annulus in $\bbr^p$. 
Furthermore, the riemannian measure on $B\times A$
(for the metric induced from $Gr\times \rn$)
is smoothly equivalent to the product measure. Hence, it suffices to consider the cartesian case.
For simplicity we drop the tildes and   rewrite the condition (\GGG.8) as
$$
  u^j(w,x) \ \to\    U(w,x) \qquad{\rm in}\ \ L^1(B\times A).
\eqno{(\GGG.8)'}
$$
and  rewrite the condition (\GGG.9) as
$$
   u^j(w,x) \ \to\       V(w,x) \qquad{\rm in}\ \ L^1(A) \qquad{\rm for\ all\ nonpolar\ \ } w\in B.
\eqno{(\GGG.9)'}
$$
  By (\GGG.8)$'$ we have that
$$
|  u^j(w,x) -   U(w,x)| \ \to\  0 \qquad{\rm in}\ \ L^1(B\times A)
\eqno{(\GGG.8)''}
$$
and by (\GGG.9)$'$ we see that for all non-polar  $w\in B$, 
$$
|  u^j(w,x) -   U(w,x)| \ \to\  |   V(w,x) -     U(w,x)|  \qquad{\rm in}\ \ L^1( A).
\eqno{(\GGG.9)''}
$$
Now by the Fubini Theorem, the function
$$
 I^j(w) \ \equiv  \ \int_A |u^j(w,x) -U(w,x)|\, dx 
$$
in integrable on $B$, and 
$$
\int_B I^j(w)\, dw\ =\ \int_{B\times A}  |u^j(w,x) -U(w,x)|\, dw\, dx.
\eqno{(\GGG.11)}
$$
Moreover,  by (\GGG.9)$''$ we know that 
$$
I^j(w) \ \ {\rm converges \ pointwise\ to }\ \ \int_A |  V(w,x) -   U(w,x)|  \, dx \ \ {\rm on}\ \ B.
\eqno{(\GGG.12)}
$$
By Fatou's Lemma and  (\GGG.12), (\GGG.11) and (\GGG.8)$''$ we have
$$
\int_{B \times A}  |V(w,x) -U(w,x)| \, dw\, dx \ =\ \int_B \lim I^j(w) \, dw \ \leq\ \underline{\lim} \int_b I^j(w)\,dw \ =\ 0.
$$
Thus  $V=U$ a.e. on $B\times A$. Furthermore, 
 for almost all $w\in B$,  the restrictions satisfy $V\bigr|_{w\times A}= U\bigr|_{w\times A}$ a.e. on $A$.
  Since these restrictions are both $\D$-subharmonic on $A$, we conclude equality everywhere
  on $A$.  This establishes (\GGG.10)
   and  completes the proof of Lemma \GGG.2.  \qed\medskip
  
  Finally we prove (\GGG.6).
 Let   $\GG^*\ss\GG$ denote the set of non-polar $\GG$-planes $W$ 
 for which 
 $$
 U\bigr|_W\  = \   \den(W) K_p
 \eqno{(\GGG.13)}
 $$
This set has full measure in $\GG$ by (\GGG.10).

Consider  a general non-polar $p$-plane $W\in \GG$.  
  Let $S_W = W \cap S^{n-1}$ denote the unit sphere in $W$.
Since   $g \equiv U\bigr|\snm$ is upper semi-continuous on $S_W$,
 it assumes its  maximum at a point $x\in S_W$.
Now since
$$
\wt U(W,x)  \ =\  {\rm ess} \!\!\! \!\!\! \limsup_{(W',x')\to (W,x)}  \wt U(W',x') 
$$
and $\GG^*$ has full measure,

\smallskip
\centerline
{ There exists a sequence $(W_j, x_j) \in Y$, with $|x_j|=1$ and $W_j\in\GG^*$}
 \centerline{
such that   $(W_j,x_j)\to (W,x)$ and $\wt U (W_j,x_j) \to \wt U (W,x)$}
\smallskip
\noindent
Choose another unit vector $y\in S_W$. Since $W_j\to W$, we have 
$S_{W_j}\to S_W$ and there exists a 
sequence of unit vectors $y_j\in W_j$ such that $y_j\to y$.  
Since $W_j \in \GG^*$, we have $g(x_j) = g(y_j) = -\den(u\bigr|_{W_j})$.
Thus, using the upper semi-continuity of $g$ we have
$$
g(x) \ =\  \lim_j g(x_j)\ =\ \lim_j g(y_j) \  \leq \ g(y),
$$
and $g(y)\leq g(x)$ since $g(x)$ is the maximum value of $g$ on $S_W$.
We have proved that $g$ is constant on $S_W$.
By Corollary \GGG.3  and the definition of density, we now conclude (\GGG.13) for
our general non-polar plane $W$. Finally note that $U\bigr|_W \equiv -\infty$
if $W$ is polar for $u$ at 0.

We have now proved Theorem \GGG.1 for $p\geq3$ under the assumption that $\GG$ is a smooth
submanifold. For a general $\GG$, choose any point $W\in\GG$ and consider the 
$G$-orbit $\GG_0 \equiv G\cdot W \ss \GG$.  Now $\GG_0$ is a compact smooth
submanifold of $G(p,\rn)$, and since $F\ss F_0\equiv F(\GG_0)$, we see that any
$F$-subharmonic function is $F_0$-subharmonic. Hence the result for smooth 
$\GG$ implies the result in general.

We now address the case $p=2$.  The proof given here follows that of Kiselman 
[K] and is easier than the one given above for the cases $p\geq3$.
Our first observation is that by the first equality in (\FF.8) of Part I, we have
$$
U(x) \ \leq\ \den^M(u) \log\,|x| \qquad x\in\rn.
$$
Suppose now that $W\in\GG$ is non-polar for $u$ at 0. Using either of the complex structures
induced on the 2-plane $W$ by the inner product, we have that if $x\in W$ and $\l\in \bbc$,
then $U(\l x) \leq  \den^M(u) \log\,|\lambda x|$.  Hence,
$$
U(\l x) - \den^M(u) \log\,|\lambda | \ \leq\  \den^M(u) \log\,|x|.
$$
For $x\in W-\{0\}$, the function $V(\l) \equiv U(\l x) - \den^M(u) \log\,|\lambda |$
is bounded above and $\D$-subharmonic on $\bbc$ (since $U\bigr|_W$ is $\D$-subharmonic).
By Liouville's Theorem $V(\l)$ is constant  equal to $V(1)=U(x)$. Hence, $V(\l)=V(1)$
says that 
$$
U(\l x) \ =\  \den^M(u) \log\,|\lambda |  + U(x) \qquad\forall\ x\in W, \ \l\in\bbc.
$$
Setting $y=re^{i\theta} \in W$, $x=e^{i\theta} $, this gives the desired result:
$$
U(y) \ =\  \den^M(u) \log\,|y |  + U\left({y\over |y|}\right) \qquad\forall\ y\in W.
$$
Now  the first equality in  (\FF.8) of Part I shows that
$
\sup_{|x|=1} U(x) \ =\ \den^M(u) \log \,1 \ =\ 0.
$
\qed

%\Note{} We know that $-g\bigr|_{W\cap S^{n-1} }$ in  Part (b) of Theorem \GGG.1 is equal to 
%the constant $\den_W$ coming from Proposition \FF.6, for almost all $W$.  We have yet to prove equality %everywhere.
%This loose end may not be important in the end.

%%%%%%%%%%%%%%%%%%%%%%%%%%%%%%%%%%%%%%%%%%%%%%%%
%%%%%%%%%%%%%%%%%%%%%%%%%%%%%%%%%%%%%%%%%%%%%%%%
%%%%%%%%%%%%%%%%%%%%%%%%%%%%%%%%%%%%%%%%%%%%%%%%
%%%%%%%%%%%%%%%%%%%%%%%%%%%%%%%%%%%%%%%%%%%%%%%%
%%%%%%%%%%%%%%%%%%%%%%%%%%%%%%%%%%%%%%%%%%%%%%%%

\vfill\eject

\centerline{\headfont \KK. The  Strong Uniqueness  Theorems.}
\medskip 

To begin we introduce the following concept.

 \Def{\KK.1}  We say that $\GG\ss G(p,\rn)$ has the {\bf transitivity property} if 
 for any two  vectors $x,y\in\rn$ there exist $W_1,...,W_k\in\GG$
with $x\in W_1, y \in W_k$ and dim$(W_i\cap W_{i+1}) >0$ for all $i=1,...,k-1$. 
\medskip

Note that if any two points $x,y\in\rn$ are contained in $W$ for some $W\in\GG$, then, of course, 
$\GG$ has the transitivity property.

%We assume $\GG$ satisfies the standing assumptions.

\Theorem {\KK.2. (Strong Uniqueness I)}  
{\sl  Assume $\GG$ has the transitivity property.  Then  strong uniqueness of tangents holds
for all $\GG$-plurisubharmonic functions.}
\pf
Let $U$  be a tangent at 0 to a $\GG$-plurisubharmonic function $u$, and suppose $p\neq2$.
 By  Theorem  \GGG.1 we know that for every $W\in \GG$,
 $U(x) = -{\den(W)\over |x|^{p-2}}$    $\forall\,x\in W$.  Thus if $W,W'\in\GG$ satisfy
 $\dim(W\cap W')\geq1$, then $\den(W)=\den(W')$. Hence, by the transitivity property, 
 $\den(W)$ is constant on $\GG$. Clearly that constant is $\den(u,0)$.
 When $p=2$ the argument is similar.
 \qed
 
\vskip .3in

\centerline{\bf Some Examples.}

\medskip
One can establish the transitivity property for the following sets $\GG$, 
and therefore one has strong uniqueness of tangents for the corresponding
$\GG$-plurisubharmonic functions.
 
 \medskip
 
 (a)  \ \ $\GG = G(p,\rn)$ ($p$-plurisubharmonic functions) for $p>1$.

 \medskip
 
 (b)  \ \ $\GG = G^\bbc(k,\bbc^n)$ (complex $k$-plurisubharmonic functions) for $k>1$ ($p=2k$).

 \medskip
 
 (c)  \ \ $\GG = G^\bbh(k,\bbh^n)$ (quaternionic $k$-plurisubharmonic functions) for $k>1$  ($p=4k$).

 \medskip
 
 (d)  \ \ $\GG = $ ASSOC (Associative subharmonic functions in $\bbr^7$) ($p=3$).

 \medskip
 
 (e)  \ \ $\GG = $ COASSOC (Coassociative subharmonic functions in $\bbr^7$)  ($p=4$).

 \medskip
 
 (f)  \ \ $\GG =$ CAYLEY (Cayley subharmonic functions in $\bbr^8$)  ($p=4$).

 \medskip
 
 (g)  \ \ $\GG = $ LAG (Lagrangian subharmonic functions in $\bbc^n$)  ($p=n$).

 \medskip
 
 (h)  \ \ $\GG = \ISO$  ($p$-isotropic subharmonic functions in $\bbc^n$).

\Note{\KK.3} In the three cases: $G(1,\rn)$ (i.e.,  $F=\cp$), 
$G(1,\bbc^n)$ (i.e.,  $F=\cp^\bbc$),  and $G(1,\bbh^n)$ (i.e.,  $F=\cp^\bbh$), 
strong uniqueness fails. In Section \LL \
the possible tangents in these cases are completely characterized. In the convex case uniqueness of tangents holds, which of course is classical.  In the complex
case, uniqueness fails. This is due to Kiselman [K]. 

\medskip
Strong uniqueness in cases (a), (b) and (c) above also follows from Theorem  \JJJ.1
in Part I.  However, the others do not.

\vfill\eject
\Theorem{\KK.3. (The Transitivity Theorem)} 
{\sl   Fix $p\geq 2$ and $n\geq 3$.

\smallskip

(a) Every compact SU$(n)$-invariant subset 
$
\GG \ss G^\bbr(p, \bbc^n)
$
except $\cp^\bbc$ has the transitivity property.
\smallskip

(b) Every compact Sp$(n)\cdot {\rm Sp}(1)$-invariant subset 
$
\GG \ss G^\bbr(p, \bbh^n)
$
with three exceptions has the transitivity property. 
 The exceptions are the sets of real $p$-planes which lie in 
 a quaternion line for $p=2,3,4$.  When $p=4$ this is $\cp^\bbh$

(c) For $p\geq5$, every compact Sp$(n)$-invariant subset 
$
\GG \ss G^\bbr(p, \bbh^n)
$
 has the transitivity property.
}

\medskip
\noindent
{\bf Proof of (a).}
When $p\geq 3$ there is a simple argument, which we give first.
Let $W\ss \bbc^n$ be a real 3-plane and consider the orbit
$
{\rm SU}(n)\cdot W \ss G^\bbr(3, \bbc^n).
$

\Lemma {\KK.4}
{\sl
Given any unit vector $x \in \bbc^n$ and any unit vector  $e\perp \span\{ x, Jx\}$, there exists 
$V\in {\rm SU}(n)\cdot W$ with $x,e\in V$.
}
\pf  Clearly there exists $W' \in {\rm SU}(n)\cdot W$ with $x\in W'$.  Let $H=\{x, Jx\}^\perp$. Then
$\dim_{\bbr}(W' \cap H)\geq1$ so there exists a unit vector $e_0\in W' \cap H$.
Thus $W' = \span\{x, v, e_0\}$  for a unit vector $v\perp x, e_0$.

%We assume $n\geq 3$ since SU(2) acts transitively on real hyperplanes in $\bbc^2$.
Choose $g\in  {\rm SU}(n)$ such that $gx=x$ and $g(e_0)=e$.  This is possible
since ${\rm SU}(n-1)\equiv \{g: gx=x\}$ acts transitively on the unit sphere in $\bbc^{n-1}$
for $n>2$.
Set $V = g(W') = \span\{x, gv, e\}$. \qed

\Cor {\KK.5}
{\sl  
For any real 3-dimensional subspace $W\ss\bbc^n$, the set 
${\rm SU}(n)\cdot W$ has the transitivity property.
Consequently, any compact ${\rm U}(n)$-invariant subset 
$
\GG\ \ss\ G^\bbr(p,\bbc^n)
$
for $p\geq3$ has the transitivity property.
}
\pf
%As above, we can assume $n\geq3$ since SU(2) acts transitively on real
%3-planes in $\bbc^2$, and any two 3-plane have a common line.
Given non-zero vectors $x,y\in \bbc^n$, choose a unit vector $e$ with
$e\perp x,Jx,y,Jy$. By Lemma \KK.4 there exist $W_x, W_y\in {\rm SU}(n)\cdot W$
with $x,e\in W_x$ and $y,e\in W_y$. Thus, ${\rm SU}(n)\cdot W$ has the transitivity property.
The second assertion follows immediately.\qed

\medskip
This leaves the case where $p=2$.  

\Lemma {\KK.6}
{\sl Given a real 2-plane $W\ss\bbc^n$ and any  (real)
orthonormal basis $x,v$ of $W$, the number $|\bra {Jx}v| \equiv\cos\theta$ 
is a complete   invariant of the orbit }
$$
\GG \ \equiv\ {\rm SU}(n)\cdot W\ \ss\ G^\bbr(2, \bbc^n).
$$
\pf Suppose $y, w$ are orthonormal  with $|\bra {Jy}w| = \cos\theta$.
There exists $g\in {\rm SU}(n)$ with $gx=y$, so we may assume $y=x$.  
By changing the sign of (say) $v$ if necessary, we may assume 
$\bra {Jx}v= \bra {Jx}w$.  Now
$v=\bra {v}{Jx} Jx + v_0$ and $w=\bra {w}{Jx} Jx + w_0$, where $v_0$ and $w_0$ 
are orthogonal to $x,Jx$. Now since $n\geq3$, there exists $g'\in  {\rm SU}(n-1)\equiv \{g: gx=x\}$, as above,  so that $g'(v_0)=w_0$, and therefore $g'(v)=w$. \qed

\medskip

Part (a) for $p=2$ is  a consequence of the following.

\Prop {\KK.7}
{\sl
If $W\in G^\bbr(2,\bbc^n)$ is not a complex line, then the orbit ${\rm U}(n)\cdot W$ has the transitivity property.
}
\pf
  Fix a unit vector $x\in \bbc^n$
and consider the set 
$$
B_x \ \equiv \{v\in \bbc^n : |v|=1 \ \ {\rm and}\ \ \bra{Jx} v =   \cos\theta\}
$$
By assumption the invariant $\cos \theta \neq1$.
Hence this set is a geodesic ball in $S^{2n-1}$ of intrinsic radius
$0<\theta<\pi$ about the  point $Jx$. 

Now suppose $y\in \bbc^n$ is another unit vector with the property that
$
\partial B_{x} \cap \partial B_{y} \ \neq\ \emptyset,
$
and choose $v\in \partial B_{x} \cap \partial B_{y}$. Then $x\in \span\{x,v\}$
and $y\in\span\{y,v\}$ and by Lemma \KK.6 both $\span\{x,v\}$ and $ \span\{y,v\}$
lie in $\GG$.  

In the event that  
$
\partial B_{x} \cap \partial B_{y} = \emptyset,
$
 we can find a sequence of points $x=x_0, x_1, x_2,...,x_N=y$
such that 
$$
\partial B_{x_{k-1}} \cap \partial B_{x_{k}} \ \neq\ \emptyset\qquad {\rm for}\ \ k=1,...,N.
$$
This completes the proof of (a)

\medskip
\noindent
{\bf Proof of (c).}  This closely follows the arguments given in Lemma \KK.4 and Corollary \KK.5,
and is omitted.

\medskip
\noindent
{\bf Proof of (b).}
Let $W\ss \bbh^n$ be a real 2-plane, and choose an orthonormal basis $\{x,v\}$ of $W$.
Let $\pi^\perp_x$  denote orthogonal projection onto $(\bbh x)^\perp\ss \bbh^n$.

\Lemma{\KK.8} {\sl
The norm $|\pi^\perp_x(v)|^2$ is independent of the choice of orthonormal basis
$\{x,v\}$ for $W$, and it is a complete invariant for the action of ${\rm Sp}(n)\cdot {\rm Sp}(1)$
acting on the Grassmannian $G^\bbr(2, \bbh^n)$.
}
\pf Let $e_0=1, e_1,e_2,e_3$ be an orthonormal basis of $\bbh$. Then 
$|\pi^\perp_x(v)|^2 = 1 -  {\bra x v}^2-\sum_{j=1}^3 { \bra v {e_jx}}^2 = 1 -\sum_{j=1}^3 { \bra v {e_jx}}^2$. Now let
$x'=x \cos \theta +v \sin\theta$ and $v'= -x \sin \theta +v \cos\theta$  be another orthonormal
basis of $W$. Using the fact that $\bra {v} {e_j x} = -\bra {e_jv}  x$,  
one computes that $\bra {v'} {e_j x'} = \bra v {e_jx} \cos^2\theta -\bra  {e_jv} x \sin^2\theta =
\bra {v} {e_j x}$ for $j=1,2,3$. This proves the independence of the choice of orthonormal basis.

Now suppose we have 2-planes with o.n. bases $W=\span\{x,v\}$ and $W'=\span\{x',v'\}$.
Then there exists $g\in{\rm Sp}(n)$ with $g(x)=x'$, so we may assume that $x=x'$.
Let $\pi_x$ denote orthogonal projection onto the quaternion line $\bbh x$.
The subgroup of ${\rm Sp}(n)\cdot {\rm Sp}(n)$ which maps $\bbh x$ to itself is 
transitive on all real 2-planes in $\bbh x$, in fact  it contains an SO(4)-subgroup acting 
standardly on $\bbh x = \bbr^4$.
%Thus by applying such a transformation we have $\span\{x, \pi_xv\} = \span\{x, \pi_xv'\}$.
Thus there is an element in this subgroup which fixes $x$ and maps $\pi_xv$ to $ \pi_xv'$
(since they are orthogonal to $x$ and have the same length). 
Now since Sp$(n-1) \equiv \{g\in {\rm Sp}(n) : g(x)=x\}$ acts transitively on the unit sphere
in $(\bbh x)^\perp$, it contains an element which maps $\pi_x^\perp v$ to $ \pi_x^\perp v'$.\qed

\Prop{\KK.9} {\sl
Let $W$ be a real 2-plane in $\bbh^n$ which is not contained in a quaternion line.
Then the orbit $\GG\equiv {\rm Sp}(n) \cdot {\rm Sp}(1) W$  has the transitivity property.
}
\pf
Let $W=\span\{x_0, v_0\}$ as above. 
By assumption the  invariant 
$\sin^2\theta \equiv |\pi_{x_0}^\perp v_0|^2\neq0$. 
Fix a unit vector $x\in\bbh^n$.
By Lemma \KK.8  the set of 2-planes in  $\GG$ which contain $x$ is 
$$
\Sigma_x \ \equiv\ \{ W\in\GG : x\in W\} \ =\  \{\span\{x,v\} : |\pi_{x}^\perp v|^2 = \sin^2\theta\}
\ \cong\ \{v\in S^{4n-1} : |\pi_{x}^\perp v|^2 = \sin^2\theta\}.
$$
This is  the real hypersurface of points in $S^{4n-1}$ at constant distance $\theta$ from the 
geodesic 2-sphere $S^2_x\equiv \{e\cdot x : e\in  {\rm Im} \bbh \ {\rm and}\  |e|=1\}$.
Now it is straightforward to see that $\Sigma_x\cap \Sigma_y\neq \emptyset$ for all 
$y$ sufficiently close to $x$.  By homogeneity the measure of closeness is independent of
$x$. The transitivity property follows.\qed

\medskip
Assertion (b) now follows, and the proof of Theorem \KK.3 is complete.
\qed\medskip

Theorem \KK.3 implies that nearly every SU$(n)$- or Sp$(n)$-invariant set 
$\GG$ has the transitivity property. Among the geometrically interesting examples
are the  sets of Lagrangian and, more generally, isotropic  planes in $\bbc^n$
(see Example A.4). Here are   further examples.

\Ex{\KK.10. (Cauchy-Riemann Sets)}
Fix integers $1\leq m$ and $p> 2m$, and define
$$
\GG \ =\ \{W\in G^{\bbr}(p, \bbc^n) : \dim_{\bbc}(W\cap JW)\geq m\}
$$
Closely  related is the set
$$
\GG_0 \ =\ \{V\oplus L \in G^{\bbr}(2m+\ell, \bbc^n) : V=JV, \ L\perp JL \ {\rm and}\ \dim_\bbc V=m \}
$$
Notice that $\GG_0$-submanifolds have constant CR-rank $\equiv m$, and 
$\GG$-submanifolds have CR-rank  $\geq m$.

\Ex{\KK.11. (Quaternionic Isotropic and Cauchy-Riemann Sets)}
In $\bbh^n$ we have the sets of isotropic $p$-planes ($p\leq n$)
$$
{\rm ISO}_p^\bbh \ \equiv\ \{W\in G^{\bbr}(p, \bbh^n) : Iw, Jw, Kw \perp W \ \forall\, w\in W\}
$$
There are also quaternionic analogues of the Cauchy -Riemann sets given in Example \KK.10.

\Theorem{\KK.12. (Strong Uniqueness II)} 
{\sl   Fix $p\geq 2$ and $n\geq 3$.
Then strong uniqueness of tangents  to  $\GG$-plurisubharmonic functions
holds for:
\smallskip

(a) Every compact SU$(n)$-invariant subset 
$
\GG \ss G^\bbr(p, \bbc^n)
$
except $\cp^\bbc$,
\smallskip

(b) Every compact Sp$(n)\cdot {\rm Sp}(1)$-invariant subset 
$
\GG \ss G^\bbr(p, \bbh^n)
$
with three exceptions, namely 
 the sets of real $p$-planes which lie in 
 a quaternion line for $p=2,3,4$  (when $p=4$ this is $\cp^\bbh$),

(c) For $p\geq5$, every compact Sp$(n)$-invariant subset 
$
\GG \ss G^\bbr(p, \bbh^n).
$
 
}
\pf This is an immediate consequence of Theorems \KK.2 and \KK.3 above.\qed

%%%%%%%%%%%%%%%%%%%%%%%%%%%%%%%%%%%%%%%%%%%%%%%%
%%%%%%%%%%%%%%%%%%%%%%%%%%%%%%%%%%%%%%%%%%%%%%%%
%%%%%%%%%%%%%%%%%%%%%%%%%%%%%%%%%%%%%%%%%%%%%%%%
%%%%%%%%%%%%%%%%%%%%%%%%%%%%%%%%%%%%%%%%%%%%%%%%
%%%%%%%%%%%%%%%%%%%%%%%%%%%%%%%%%%%%%%%%%%%%%%%%

\vfill\eject

\centerline{\headfont \EE.\ Homogeneous F-Subharmonics.}
 \medskip 
 
 We begin by computing the  formula for the second derivative $D_x^2 u$  of a 
 function $u$, which is homogeneous of degree $m$, in terms of its restriction 
 $g\equiv u\bigr|_{S^{n-1}}$ to the unit sphere $S^{n-1}\ss\rn$.
 For our application it is useful to replace $m$ by $p\equiv -m+2$ (or $m=-(p-2)$) so that:
 $$
 u(x) \ =\ {1\over |x|^{p-2}} \ g\left(  {x\over|x|}  \right).
 \eqno{(\EE.1)}
$$

\noindent
{\bf Remark \EE.1.  ($p=2$).}  In the special case $p=2$
the natural extension of $g$ is given by $u(x) =\den \log|x| +  g\left(  {x\over|x|}  \right)$.
This choice is consistent with the Riesz kernels and with classical pluripotential theory.
The formulas computed below, when $p=2$, only apply to the special case $\den=0$.
However, they also apply directly to give the corresponding formulas in the general case.
This is discussed in Remark \LL.3.

\medskip

 Let $\Hess_\s  g$ denote the riemannian hessian of $g$ at a point  $\s=x/|x| \in S^{n-2}$.
 Using the orthogonal  decomposition 
 $$
 \rn \ =\  T_\s(\rn) \ = \ T_\s(S^{n-1}) \oplus N_\s(S^{n-1})
 \eqno{(\EE.2)}
$$
  the quadratic form $\Hess_\s g$ on $T_\s(\rn)$
 can be considered to be a quadratic form on $\rn$
 (whose null space contains $N_\s(S^{n-1})$).
 Also, the tangential derivative
 $D_\s g=dg$ at $\s$ can be considered a vector in $\rn$.  Then
 $$
|x|^p D_x^2 u \   =\  \Hess_\s g  - (p-2) g P_{x^\perp} - (p-1)  \left(\s \circ D_\s g \right) 
+ (p-2)(p-1) g P_x    
\eqno{(\EE.3)}
$$
where on the right hand side,  $J^2g \equiv (g, Dg, \Hess\, g)$, the riemannian 2-jet of $g\in C^2(S^{n-1})$,
is evaluated at the point $\s= x/|x| \in S^{n-1}$.

\medskip
\noindent
{\bf Proof of (\EE.3).}  
One computes that   $D ( {1\over |x|^{p-2}}) = - {p-2\over |x|^{p-1} } \,{x\over|x|}$, and from Lemma \AA.1,
we have that 
$$
|x|^p  D^2 \left ( {1\over |x|^{p-2}}\right ) = -(p-2) P_{[x]^\perp} +(p-2)(p-1) P_{[x]}.
\eqno{(\EE.4)}
$$
Define $\wt g(x) = g(x/|x|)$ for $x\in \rn-\{0\}$.  Then direct calculation shows that 
$D \wt g(x) \ \cong \ (1/|x|) D_\s g$ 
and that
$$
D^2_x \wt g \ \cong \ {1\over |x|^2} \left (\Hess\, g - {\s}\circ Dg\right)
\eqno{(\EE.5)}
$$
These formulas for the first and second derivatives of  the functions
$1/|x|^{p-2}$ and $\wt g$ yield the formula (\EE.3) for the second derivative of the product $u$.\qed

\medskip

Define $\Phi(J^2_\s g) \in \Symn$ to be the RHS of (\EE.3).  That is,
$$
\Phi(J^2_\s g) \ \equiv  \ \Hess_\s g - (p-2) g(\s) P_{\s^\perp} - 
(p-1) ( \s \circ D_\s g) + (p-2)(p-1) g(\s) P_\s.
\eqno{(\EE.6)}
$$
Then (\EE.3) says that 
$$
|x|^p D^2_xu \ =\ \Phi(J^2_\s g) \quad{\rm with}\ \ \s\equiv {x\over |x|}.
\eqno{(\EE.3)'}
$$

In terms of the $2\times 2$-blocking induced on $\Symn$ by the decomposition
(\EE.2) 
$$
\Phi(J^2_\s g) \  = \ \left(
\matrix
{
\Hess\, g   -(p-2) g I &    -(p-1) Dg  \cr
\ &\ \cr
 -(p-1) Dg^t  &   (p-2)(p-1) g
}
\right)
\eqno{(\EE.6)'}
$$
with RHS evaluated at $\s$ (and with $Dg$  written as a column vector.)

The formula (\EE.3)$'$ has been proved for $u$ and $g$ related by (\EE.1) and of
class $C^2$.  This immediately implies the following.

\Prop{\EE.2}  {\sl For a cone subequation $F$ and $u(x) = {1\over |x|^{p-2}} g({x\over |x|})$ of class $C^2$,
\medskip
\centerline
{
 $u$ is $F$-subharmonic on $\rn-\{0\} 
\qquad\iff\qquad
\Phi(J_x(g))\in  F  \ \ \ \forall\, |x|=1$, \qquad\qquad {\rm and}
}
\medskip
\centerline
{
\quad $u$ is $F$-harmonic on $\rn-\{0\} 
\qquad\iff\qquad
 \Phi(J_x(g))   \in\partial F  \ \ \ \forall\, |x|=1$,    \qquad\qquad
}
}
\medskip

We wish to extend this proposition to include upper semi-continuous functions
$u$ and $g$.  Note that with $u$ and $g$ related by (\EE.1), $u$ is upper semi-continuous on 
$\rn -\{0\}$ if and only if   $g$ is upper semi-continuous on $S^{n-1}$.

\Lemma{\EE.3} {\sl Given a subset $F\ss \Symn$, consider the subset 
$$
F_{S^{n-1}}\ \equiv\ \Phi^{-1}(F)
\eqno{(\EE.7)}
$$
of the 2-jet bundle $J^2(S^{n-1})$.

\def\imp{\qquad\Rightarrow\qquad}

\medskip

(1)\ \ $F$ closed $ \imp F_{S^{n-1}}$  closed.

\smallskip

(2)\ \ $F$ satisfies (P) $\imp  F_{S^{n-1}}$  satisfies (P).

\smallskip

(3) \ \ $F$ is a cone $\imp F_{S^{n-1}}$ is a cone bundle,

\smallskip

(4)\ \ $F$ is a convex cone $\imp F_{S^{n-1}}$ is a convex cone bundle.

\smallskip

(5) \ \ For any subgroup $H\ss {\rm O}(n)$
\smallskip
\centerline{
$F$ is H-invariant  $\imp F_{S^{n-1}}$ is H-invariant }

\smallskip

(6)\ \  If $F$ is a cone subequation, i.e., (1), (2) and (3) are true,  then the dual
$$
\wt{F_{S^{n-1}}} \ =\  \left( \wt F\right)_{S^{n-1}}.
$$

\smallskip

(7)\ \ Suppose $F$ is a cone subequation with Riesz charactersitic $p$.
\smallskip

\indent\indent (a)\ \ If $p>2$, then
\smallskip
\centerline
{
$F_{S^{n-1}}$ satisfies (N) $ \qquad\iff\qquad $  $F$ is $\cp_p$-monotone.
}\smallskip

\indent\indent (b)\ \ If $1\leq p<2$, then (N) fails for $F_{S^{n-1}}$.
}
\medskip

Before proving this lemma we state the main result. But first consider the following.

\Ex{\EE.4} Let $F\equiv \{A : \tr \,A\geq0\} =\D$ be the standard Laplacian on $\rn$.
Then by (\EE.6)$'$ $F_{S^{n-1}}$  is the linear subequation $Lg\geq0$ on $S^{n-1}$, where
$$
Lg \ \equiv \tr \, \Phi(J(g)) \ =\ \D_{S^{n-1} } g -  (n-p)(p-2)g.
$$
Note that $L$ satisfies (N) if  $2<p\leq n$.

\medskip

Now the extension of Proposition \EE.2 to include u.s.c. functions $u$ and $g$ can be stated as follows.
\vfill\eject

\Theorem{\EE.5. ($p\neq 2$)} {\sl
Suppose that $F\ss\Symn$ is a cone subequation.
If $u$ and $g$ are upper semi-continuous functions related by (\EE.1), then

\medskip
\centerline
{
 $u$ is $F$-subharmonic on $\rn-\{0\}
\qquad\iff\qquad
\ \  g$ is $F_{S^{n-1}}$-subharmonic on $S^{n-1}$,\ \ and
}
\medskip
\centerline
{
 \qquad $u$ is $F$-harmonic on $\rn-\{0\} 
\qquad\iff\qquad
g$ is $F_{S^{n-1}}$-harmonic on $S^{n-1}$ \qquad\qquad\quad
}
\medskip
\centerline
{
$\qquad\qquad\qquad\qquad\qquad\qquad  \qquad\quad     \      \iff\qquad
g$ is $F_{S^{n-1}}$-subharmonic   on $S^{n-1}$ and 
}
\smallskip
\centerline
{
$\qquad\qquad\qquad\qquad\qquad\qquad\qquad\qquad\qquad \
 -g$ is $\ft_{S^{n-1}}$-subharmonic on $S^{n-1}$.
}
}
\medskip

In the  applications typically  $u$ is $F$-subharmonic across 0. 
This imposes an additional condition on $g$.

\Prop{\EE.6}  {\sl
Suppose $F$ has Riesz characteristic $1 \leq p<\infty$, $p\neq 2$, 
and that $u(x) \equiv {1\over |x|^{p-2}} g({x \over |x|})$ is $F$-subharmonic on $\rn$ across 0.
Let  $\den \equiv \den^M(u,0)$  be the density of $u$ at 0.
Then:
\medskip

(1)\ \ $g$ is $F_{S^{n-1}}$-subharmonic on $S^{n-1}$, and 
\medskip

(2)\ \ if $2< p<\infty$, then $\sup_{S^{n-1}}g \ =\ -\den \ \leq\ 0$,  while
\medskip

(3)\ \  if   $1 \leq p< 2$, then $\sup_{S^{n-1}}g \ =\ \den \ \geq\ 0$.
}
\pf  Note that $K(1)= -1$ for $2 < p< \infty$, while 
$K(1)= 1$ for $1 \leq p< 2$.  Now Proposition \BB.6 can be used to compute $\den^M(u,0)$.\qed

\medskip
\noindent
{\bf Remark \EE.7.}  If $F$ is $\cp_p$-monotone and $2<p$, then a converse is true, since (2) $\Rightarrow$ $u$ is locally bounded above across 0, in which case the singularity at 0
 is removable by results in [\HLRS]. Thus (1) and (2) imply that $u$ is $F$-subharmonic on $\rn$.

\medskip
\noindent
{\bf Proof of Lemma \EE.3.}
  Formula (\EE.6)$'$ shows that 
$$
\eqalign
{
\Phi : J^2(&S^{n-1}) \ \arr\   S^{n-1}\times \Symn \equiv \Symn \bigr|_{S^{n-1}}  \cr
&{\rm is \ an \  O(n)-equivariant  \  bundle\  map,}
}
\eqno{(\EE.8)}
$$
  and, in fact when
$(p-1)(p-2)\neq 0$,  it is  a bundle isomorphism.  
From this the implications   (1), (3), (4) and (5) are obvious.
To prove (2) note that with 
$$
J_\s(S^{n-1}) \cong \bbr\times T_\s^* S^{n-1}\times\Sym(T_\s^* S^{n-1}),
$$
 we have
$$
\Phi(0,0,P) \ =\ \left(\matrix{ P & 0 \cr 0 & 0\cr}\right) \ \in\ \cp \qquad
{\rm if}\ \ P\geq0.
$$
 To prove (6) note that 
(\EE.8) implies that $\Phi^{-1}(\Int F) = \Int \Phi^{-1}(F)$, and that the fibres of
$\Int \Phi^{-1}(F)$ are the fibres of  $\Phi^{-1}(\Int F)$.
For (7) note that by (\EE.6)
$$
\Phi(J+(-r,0,0)) = \Phi(J) + r(p-2) (P_{e^\perp} - (p-1) P_e).
$$
Hence, for $p>2$, $F_{S^{n-1}}$ satisfies (N) $\iff$ $F$ is $\cp_p$-monotone,
while for $1<p<2$, $F_{S^{n-1}}$ doesn't satisfy (N).\qed

\medskip
\noindent
{\bf Proof of Theorem \EE.5.}  The implications $\Rightarrow$ are easy since a test function
$\psi$ for $g$ at $\s$ induces a test function $\vf(x) \equiv  {1 \over   |x|^{p-2} } \psi({x\over |x|})$ for $u$
at $\s$.

To prove the reverse implications we fix a point $x_0$ which we may
assume  to be of the form $x_0 = (\rho, 0, ..., 0)$ for $\rho>0$. We then
choose the  local coordinate $y$ on the sphere about $(1, 0,...,0)$
given by $\Psi(y) = (1,y)/|(1,y)|$ for $|y|<1$. Setting $t=r-\rho$ gives 
local coordinates $(t,y)$ about $x_0$ on $\rn$ with $(0,0)$ corresponding to $x_0$. 

Under this coordinate change a function of the form $|x|^{2-p}g({x\over |x|})$ becomes
$(\rho + t)^{2-p}\g(y)$. To complete the proof it will suffice to prove the following lemma.

\Lemma{\EE.8} {\sl 
Suppose $\vf(t,y))$ is a strict quadratic test function for the function $u(t,y) = (\rho + t)^{2-p}\g(y)$ at $(0,0)$.
Then there exists a smooth test function $\psi(y)$ for 
$\g(y)$ at $0$ in  $\bbr^{n-1}$ such that }
$$
(\rho + t)^{2-p}\g(y) \ \leq\ (\rho + t)^{2-p}\psi(y) \ \leq\ \vf(t,y)\qquad{\rm near \ } (0,0).
$$
\pf
We can assume that $\g(0)=0$.  By assumption $\vf$ is a strict test function of the form
$$
\vf(t,y) \ =\ pt + \bra q y +at^2 +2  \bra b y  t+ \bra {Cy}y.
$$
Setting $y=0$ gives $0=u(t,0) < pt + at^2$ and therefore 
$$
p\ =\ 0  \and  \ a\ >\ 0.
$$
We now have
$$
  (\rho + t)^{2-p}\g(y)  \ \equiv\ u(t,y) \ \leq \ \vf(t,y) \ =\ 
at^2 +2  \bra b y  t+ k \ \equiv \ Q_y(t)
$$
where 
$$
k \ = \ k(y) \ \equiv\ \bra q y + \bra {Cy}y.
$$
For fixed $\d>0$ small we define
$$
\psi(y)\ \equiv\ \inf_{|t|\leq \d}  {1  \over  (t+\rho)^{2-p}} Q_y(t).
$$
Then on $\{|t|\leq \d\}$ we have
\medskip

(1) \ \ $u(t,y) \ \leq \ (t+\rho)^{2-p}\psi(y)$ \qquad (because $g(y)\leq \psi(y)$), and

\medskip

(2) \ \ $(t+\rho)^{2-p}\psi(y)\ \leq\ \vf(t,y)$.

\medskip

It remains to show that $\psi(t)$ is smooth when $\d$ is taken sufficiently small.
One calculates that $t$ is a critical point of the function $t\mapsto (t+\rho)^{p-2} Q_y(t)$
in the range $t+\rho >0$
if and only if 
$$
At^2 + 2Bt + C=0
\eqno{(\EE.9)}
$$
 where
$$
A\ =\ ap,\qquad B\ =\ \rho a + (p-1) \bra b y, \qquad C\ =\ 2\rho \bra by +(p-2)k(y).
$$
When $y=0$, we have that $apt^2 +2\rho a t =0$ which happens iff
$$
t\ =\ 0\qquad{\rm or}\qquad t\ =\ -{2\rho\over p}.
$$
We choose $\d << 2\rho/p$ to rule out the second possibility.
The roots $t_1(y)$ and $t_2(y)$ of (\EE.9) with $t_1(0)=0$ are
two smooth functions of $y$ in a neighborhood of $0$.  

It remains to show that $\psi(y) = t_1(y)$.  Since $t_1(y)$ and $t_2(y)$
are the critical points of   $(t+\rho)^{p-2}Q_y(t)$, this means we must show that
$\inf_{|t|\leq\d}  (t+\rho)^{p-2} Q_y(t)$ is not assumed for $t=\pm \d$.  One checks that
this is true for $y=0$ and therefore for all sufficiently small $y$.\qed

\vskip.3in
%\vfill\eject

%%%%%%%%%%%%%%%%%%%%%%%%%%%%%%%%%%%%%%%%%%%%%%
%%%%%%%%%%%%%%%%%%%%%%%%%%%%%%%%%%%%%%%%%%%%%%
%%%%%%%%%%%%%%%%%%%%%%%%%%%%%%%%%%%%%%%%%%%%%%
%%%%%%%%%%%%%%%%%%%%%%%%%%%%%%%%%%%%%%%%%%%%%%
%%%%%%%%%%%%%%%%%%%%%%%%%%%%%%%%%%%%%%%%%%%%%%

\centerline{\headfont \LL. Tangents to Convex, C-Plurisubharmonic,  }\smallskip

\centerline {\headfont and H-Plurisubharmonic Functions.}
 \medskip

We now give a brief discussion of three  geometric cases where uniqueness 
of tangents does not hold.  These are the convex functions (where $\GG = G(1,\rn)$),
the classical complex plurisubharmonic functions  (where $\GG = G(1,\bbc^n)$), and
the quaternionic  plurisubharmonic functions  (where $\GG = G(1,\bbh^n)$). 
The results in the first case follow from classical convex analysis [R].
Those in the complex case are due to Kiselman [K]. The results in the quaternionic
case are new.

\bigskip

\centerline{\bf Tangents to Convex Functions}
\medskip

Suppose $u$ is a convex function defined in a neighborhood of 0 in $\rn$, or 
equivalently, $u$ is $\cp=F(G(1,\rn))$-subharmonic.  The Riesz characteristic
of this subequation is 1, and the appropriate homotheties are
$$
u_r(x) \ \equiv\ {1\over r} \biggl( u(rx) - u(0)\biggr), \qquad r>0.
$$
{\sl Tangents are unique.}   In fact,
$$
u_r \ \downarrow\ U\ \equiv\ \lim_{r\downarrow0} u_r \ \  \ \ {\rm uniformly\ on\ compact\ subsets\ of\ }\rn.
\eqno{(\LL.1)}
$$ 
This is easy to see geometrically.  The mappings $\Psi_r:\bbr^{n+1} \to :\bbr^{n+1}$
given by $(x,t) \mapsto {1\over r} (x, t-u(0))$ carry the epigraph of $u$ to the epigraph 
of $u_r$.  Convexity implies that for $0<r<s$, ${\rm epi}(u_r) \supset {\rm epi}(u_s)$.
The epigraphs ${\rm epi}(u_r)$ increase to ${\rm epi}(U)$, that is, the functions $u_r$
decrease to $U$. The local uniform convergence follows. 
\smallskip

{\sl Tangents are homogeneous of degree 1,  that is,}
$$
U_r(x) \ =\ U(x), \ \ {\rm i.e., }\ \ U(rx)\ =\ rU(x).
\eqno{(\LL.2)}
$$ 
This is immediate since tangents are unique.

The {\sl subdifferential} of $u$ at 0, denoted $(\partial u)(0)$,
is the set of $p\in\rn$ such that $u(x)-u(0) \geq \bra px$  for $|x|$ small.
It is easy to see that $(\partial u)(0)$ is a non-empty compact convex set. 
Now the unique tangent function $U$ to $u$ at 0 is related to the subdifferential
by
$$
U(x) \ =\ \sup_{p\in (\partial u)(0)}  \bra px
\eqno{(\LL.3)}
$$
Note that  $p\in (\partial u)(0) \iff {1\over r}(u(rx)-u(0))\geq {1\over r} \bra p{rx} = \bra px$,
and hence
$$
\eqalign
{
p\in (\partial u)(0)\quad &\iff \quad p\in (\partial u_r)(0)\quad \iff \quad p\in (\partial U)(0)  \cr
&{\rm and}\quad U(x)\ \geq\ \bra px \quad \forall \  p\in (\partial u)(0).
}
\eqno{(\LL.4)}
$$

%\vfill\eject

Finally we show that 
$$
\den^S(u,x) \ =\ 0 \iiff  u\ \ {\rm is\ differentiable\ at\ \ } x.
\eqno{(\LL.5)}
$$
(When this holds, $\den^M(u,x) = |D_xu|$.)
Both assertions in (\LL.5) remain unchanged if we subtract an affine function
from $u$. By subtracting a supporting affine function we may assume that $u\geq0$ and $u(x)=0$.
Then by (\LL.1) $U$ is differentiable at $x$ if and only if the tangent $U$ at $x$ is $\equiv 0$.
%For convenience let $x=0$.  
Now $\den^S(U,0)=\den^S(u,x)$, so if $\den^S(u,x)=0$,
then the homogeneity of $U$ implies that $\intave S U(t\s)\, d\s = \den^S(U,0) =0$ for all $t\geq0$.
However, $U\geq0$ since $u(rx)/r \downarrow U$ and $u\geq0$, and so  $U=0$.
Conversely, if $U=0$, then since $u(rx)/r$ converges uniformly to $U=0$, $u$ is differentiable at x with $D_xu=0$.

\bigskip

\centerline{\bf Homogeneous Convex Functions}
\medskip

Every convex function $U$ which is homogeneous of degree 1 is , of course, the unique tangent to itself
at 0.  By subtracting off an affine function, one can always assume that 
$$
U\geq0 \and U(0)\ =\ 0.
$$
Such functions  are classically understood.
Rewrite $U$ as $U(x)=\|x\|$. Then
$$
\|  \l x\|\ =\ \l\|x\|  \quad\forall\, \l \geq0, \ x\in \rn\  \and \|x+y\|\ \leq\ \ \|x\|+\|y\|,
\eqno{(\LL.6)}
$$
that is, $\|\bullet\|$ is a semi-norm  on $\rn$ (not necessarily balanced).  
 By (\LL.3) the unit ball $\|p\|^*\leq 1$ in the dual norm $\|\bullet\|^*$ 
is the subdifferential $(\partial U)(0)$.

 Let $U(x)$  be  a $C^2$-function which is homogeneous of degree 1, i.e., 
 $U(x) = |x| g({x\over |x|})$ where $g\equiv U\bigr|_{S^{n-1}}$.  Then
  formula (\EE.3) with $p=1$ states that 
  $$
  D^2_x U\ =\ {1\over |x|}\bigr(\Hess_e g + g(e) P_{e^\perp}\bigr)
  \qquad{\rm where\ \ } e={x\over |x|}.
  \eqno{(\LL.7)}
$$
That is, $D^2_xU$ is the pull back of the quadratic form $\Hess_e g +gI$ on the tangent space
$T_eS^{n-1}$ to the sphere using the splitting $\rn = T_eS^{n-1}\oplus \bbr\cdot e$.
 Theorem \EE.5 gives the following. 

\Prop{\LL.1}  {\sl
Let $g \in C(S^{n-1})$ be a continuous non-negative function on the sphere $S^{n-1}$, and extend
$g$ to a homogeneous function $U(x) \equiv |x|g({x\over|x|})$ of degree 1 on $\rn$.  Then}
$$
U\ \ {\rm is\ convex\ on\ }\ \rn  
\qquad\iff\qquad \Hess\, g + gI\ \geq\ 0 \ \ {\rm (in\ the\ viscosity\ sense)\ on\ \ } S^{n-1}.
$$
%The function $U$ is convex on all of $\rn$ if, in addition, $g$ satisfies $g(e) \geq \bra ep$
%on $S^{n-1}$ for some vector $p\in\rn$.}

\medskip

When $n=2$, this is the subquation $g''(\theta) + g(\theta)\geq0$ on $S^1$.
Note that the negativity condition (N)   fails.

\medskip
\noindent
{\bf Summary.}  Tangents are unique;  strong uniqueness fails; but tangents can be characterized by 
(\LL.6) or  Proposition \LL.1.

%%%%%%%%%%%%%%%%%%%%%%%%%%%%%%%%%%%%%%%%%%%%%
%%%%%%%%%%%%%%%%%%%%%%%%%%%%%%%%%%%%%%%%%%%%%
%%%%%%%%%%%%%%%%%%%%%%%%%%%%%%%%%%%%%%%%%%%%%
%%%%%%%%%%%%%%%%%%%%%%%%%%%%%%%%%%%%%%%%%%%%%
%%%%%%%%%%%%%%%%%%%%%%%%%%%%%%%%%%%%%%%%%%%%%

\vskip .3in
\centerline{\bf Tangents to Plurisubharmonic Functions in $\bbc^n$}
\medskip

In 1988  Christer Kiselman  proved that tangents to plurisubharmonic 
functions are not unique.  In fact he completely characterized the subsets of
$\lloc(\bbc^n)$ which arise as the tangent sets to psh functions. (See Theorem 4.1 in [K].)
We present  those results here.

Since the Riesz characteristic in this case is 2, the appropriate homotheties are
$$
u_r(x) \ \equiv\ u(rx) -\sup_{B_r} u, \qquad r>0.
$$
The following is (essentially) one of Kiselman's results in [K].
Let 
\medskip
\centerline{
 $\pi:\bbc^n -\{0\} \to \bbp^{n-1}_\bbc$
}
\smallskip\noindent
denote the standard  map to complex projective space, and let $\omega$ denote the
standard K\"ahler form on $\bbp^{n-1}_\bbc$ so that 
$\pi^*\omega = i \partial\dbar \log |z|$ on $\bbc^n -\{0\}$.

\Prop{\LL.2} {\sl
Suppose $U$ is a tangent to a plurisubharmonic function $u$ defined in a neighborhood
of the origin in $\bbc^n$.  Then $U$ is of the form
$$
U(x) \ =\ \den \log |x| +\pi^*g \qquad {\rm with}\quad g\in \USC(\bbp^{n-1}_\bbc)
 \eqno{(\LL.8)}
$$
where
\smallskip

(i)\ \ 
 $\den \equiv \den^M(u)$ is the (maximum) density of $u$ at 0, 
 \smallskip
(ii)\ \ 
$g$ is $\den$-quasi plurisubharmonic on $\bbp^{n-1}_{\bbc}$, that is
$$
i \partial \dbar g +\den \,\omega \ \geq\ 0,  
 \eqno{(\LL.9)}
$$

(iii)\ \ and
$$
\sup_{\bbp^{n-1}_\bbc} g \ =\ 0.
 \eqno{(\LL.10)}
$$
}
\medskip

Since only (\LL.9) is not stated in  [K], we include its straightforward proof.
Assume (\LL.8) has been established, where $g\in\USC(\bbp^{n-1}_\bbc)$.
Then 
$ i \partial\dbar U =   i \partial\dbar(\den  \log |x| +\pi^*g) = 
\pi^* ( \den\,\omega + i \partial\dbar g)$, from which one concludes that 
$$
U\ \ {\rm is\  plurisubharmonic\ on\ \ }  \bbc^n-\{0\} 
\quad\iff\quad
i \partial\dbar g +  \den\,\omega \ \geq\ 0 \ \ \  {\rm on}\ \ \bbp^{n-1}_\bbc
\eqno{(\LL.11)}
$$

\medskip\noindent
{\bf Remark \LL.3.}
This result can be deduced from the case $p=2$   in the last section (see Remark \EE.1).    If
$$
U(x) \ =\ \den \log |x|  +    \overline g\left({x\over |x|}\right) \qquad \forall\, x\in\bbc^n,
$$
then
$$
D^2_x U \ =\ {1\over |x|^2}  \left(
\matrix
{
\Hess_e  \overline g +\den I &  -D_e   \overline g\cr
-(D_e   \overline g)^t  & -\den
}
\right)
%\eqno{(\LL.9)}
$$
One can show that  the hermitian symmetric part  $(D_x^2U)_\bbc$ vanishes on $\bbc e$
and equals  $( \Hess_e \overline g)_\bbc + \den P_{(\bbc e)^\perp}$ on $(\bbc e)^\perp$
(compare the more complicated quaternionic case below). This completes a second proof.

\medskip
\noindent
{\bf Remark \LL.4.}   Note also that each $u$ is maximal on $\bbc^n-\{0\}$ since its restriction
to each complex line through the origin is $\D$-harmonic.

Proposition \LL.2 characterizes the possible tangent functions to $u$ at 0.
Kiselman's characterization of the possible tangent sets $T_0u$ can be stated as follows.

\Theorem{ \LL.5. (C. Kiselman  [K])}
{\sl  \ \ 
Suppose  that $\den\geq0$ and that  $M$ is a non-empty subset 
of the $\den$-quasi-plurisubharmonic functions on
$\bbp_\bbc^{n-1}$ with each element $g\in M$ satisfying (\LL.10).  If $M$ is closed and connected in $\lloc(\bbp_\bbc^{n-1})$, then 
  there exists a plurisubharmonic function  $u$ defined on a neighborhood
of the origin  in $\bbc^n$ such that }
$$
T_0u\ =\ M.
$$

%%%%%%%%%%%%%%%%%%%%%%%%%%%%%%%%%%%%%%%%%
%%%%%%%%%%%%%%%%%%%%%%%%%%%%%%%%%%%%%%%%%
%%%%%%%%%%%%%%%%%%%%%%%%%%%%%%%%%%%%%%%%%
%%%%%%%%%%%%%%%%%%%%%%%%%%%%%%%%%%%%%%%%%
%%%%%%%%%%%%%%%%%%%%%%%%%%%%%%%%%%%%%%%%%

%\vfill\eject
\vskip.3in

\centerline
{
\bf  Homogeneous Quaternionic Harmonics }

\medskip

The remaining series of \ST geometric cases where strong uniqueness fails
is the case of quaternionic plurisubharmonic functions ($\GG = G^\bbh(1, \bbh^n)$
in Section \KK).  Such functions have been studied by S. Alesker and M. Verbitsky
in [A$_1$], [AV], and also by the authors [\HLDD], [\HLDDR]. Note that in this case the Riesz characteristic is 4.
Let  $\pi:\bbh^n -\{0\} \to \bbp^{n-1}_\bbh$
denote the standard  map to quaternionic  projective space

\Prop{\LL.6} {\sl
Suppose $U$ is a tangent to a quaternionic plurisubharmonic function 
$u$ defined in a neighborhood of the origin in $\bbh^n$.  Then $U$ is of the form
$$
U(x) \ =\  {1\over |x|^2}  \pi^*g \qquad {\rm with}\quad g\in \USC(\bbp^{n-1}_\bbh)
 \eqno{(\LL.12)}
$$
where $g$ satisfies the subequation
$$
{\Hess}_\bbh  (g)  -2g I  \ \geq\ 0 \qquad  {\rm on}\ \ \bbp^{n-1}_\bbh.
 \eqno{(\LL.13)}
$$
Here ${\Hess}_\bbh (g)$ is the quaternionic hermitian symmetric part of 
the riemannian hessian $H=\Hess\, g$ on $\bbp^{n-1}_\bbh$, defined by 
$$
H_\bbh(v,w) \ =\ \smfrac 14  \left\{   
H(v,w)  + H(Iv, Iw) + H(Jv, Jw) + H(Kv, Kw).
\right\}
$$}

\pf
By Theorem \GGG.1 we know that $U\bigr|_{W\cap S^{4n-1}}$ is constant for every quaternion
line $W \ss \bbh^n$.  Hence, $U\bigr|_{W\cap S^{4n-1}} = \pi^* g$ for some
$g\in \USC(\bbp^{n-1}_\bbh)$. To simplify notation and to make accord with Section \EE, we shall denote
$\pi^* g$ simply by $g$.

Consider the unit sphere  $S^{4n-1} \ss \bbh^{n} =\bbr^{4n}$.  At any $x\in S^{4n-1}$ we get a decomposition
of $\bbr^{4n}$ as
$$
\bbr^{4n} \ =\ \ch_x \oplus \cv_x\oplus \bbr x
 \eqno{(\LL.14)}
$$
where $\cv_x$ is the tangent space to the fibre at $x$ of the fibration
$$
\pi : S^{4n-1} \ \arr\  \bbp_{\bbh}^{n-1}
%\eqno{(2)}
$$
and $\ch_x$ is the orthogonal complement of $\cv_x$ in $T_x S^{4n-1}$.  $\ch$ is {\sl horizontal}
for $\pi$ and $\pi_*$ maps it isometrically onto $T_{\pi x}\bbp_{\bbh}^{n-1}$. It is an $\bbh$-linear subspace of $\bbh^{n}$. Note that the radial $\bbh$-line  through $x$ satisfies
$\bbh x = \cv_x\oplus \bbr x$.

\medskip

Let $I,J,K$ be the standard basis for the imaginary quaternions.  Then we get
a trivialization of $\cv$ by the global vector fields:
$$
V_1(x) \ =\ I(x), \ \ V_2(x) \ =\ J(x), \ \ V_3(x) \ =\ K(x).
 \eqno{(\LL.15)}
$$

Now we are considering  the operator from (\EE.6)$'$
$$
L(g)  \  = \ \left(
\matrix
{
\Hess\, g   -2 g I &    -3 Dg  \cr
\ &\ \cr
 -3 Dg^t  &   6 g
}
\right)
$$
written with respect to the splitting $T_x S^{4n-1} \oplus \bbr x$. We want to compute
this for a function $g = \pi^* \wt g$ where $\wt g:\bbp_{\bbh}^{n-1} \to \bbr$. More precisely we
want to compute the quaternionic hermitian symmetric part:
$$
L_\bbh(g)(v,w) \ \equiv\ \smfrac 1 4 \{ (Lg)(v,w) + (Lg)(Iv,Iw) + (Lg)(Jv,Jw) + (Lg)(Kv,Kw) \}.
$$
Our first observation is that 
$$
\cv \ \ss\ {\rm Ker}(Dg)\and  \cv \ \ss\ {\rm Null}(\Hess g).
$$
Thus if $v\in \cv$, then we find that
$$
4 L_\bbh(g)(v,v) \ =\ -2g -2g -2g +6g \ =\ 0.
$$
Similarly, if $e=x$ we have $L_\bbh(g)(e,e) =0$, and $L_\bbh(g)(v,e) =0$.  Thus $\bbh \cdot x$
(the radial $\bbh$-line) lies in the kernel of $ L_\bbh(g)$ (cf. (\LL.17) below).

We now want to compute the spherical hessian Hess$(g)$.  Let $H$ be an invariant horizontal
vector field obtained by lifting a vector field $\wt H$ on $ \bbp_{\bbh}^{n-1}$ via $\pi$, and let
$V$ be a vertical vector field which is a real  linear combination of the $V_j$ above.
Then
$$
0 \ =\ \cl_V(H) \ =\ [V,H] \ =\ \nabla_V H - \nabla_H V
%\eqno{(5)}
$$
where $\nabla$ is the riemannian connection on the sphere.  
Observe now that
$$
 \nabla_H V \ =\  \nabla_H (Jx) \ =\ \{ \nabla^{\bbr^{4n+4}}_H (Jx)\}^{\rm Tan}
 \ =\  \{  (JH)\}^{\rm Tan} \ =\   (JH)\ \cong\ VH
\eqno{(\LL.16)}
$$
where we identify $J\cong V$ via the global identification ((\LL.15)) above:
$$
\cv \cong {\rm Im}\bbh
\eqno{(\LL.15)'}
$$

Now let $H_1, H_2$ be invariant horizontal vector fields as above.
Then
$$
\eqalign
{
(\Hess \,g) (H_1, H_2) \ &=\ H_1H_2 g - (\nabla_{H_1} H_2)g   \cr
 &=\ H_1H_2 g - (\nabla_{H_1} H_2)^{\ch}g    \cr
 & = \ \left( \Hess^{\bbp^n_\bbh} g\right) (H_1,H_2) \quad {\rm lifted\  to\  the\  sphere}
}
$$
We must now compute the $(\cv,\ch)$-component of $\Hess \,g$.
$$
\eqalign
{
(\Hess \,g) (H, V) \ &=\ HV g - (\nabla_{H} V)g   \cr
 &=\ 0-  \bra {(\nabla_{H} V)} {Dg}   \cr
 & = \ -  \bra {VH} {Dg} \ =\ -  \bra {H} {\mu^*_V Dg}
}
$$
by ((\LL.16)) above, where $\mu_V$ is the action of the imaginary quaternion $V$ on $\ch$ at
$x$.  

So with respect to the splitting ((\LL.14)) we have 
$$
L(g)  \  = \ \left(
\matrix
{
 \Hess^\bbp  g   -2 g I &   -\mu^* Dg &   -3 Dg  \cr
\mu^* Dg^t  \ &\  -2gI & 0\cr
 -3 Dg^t  &  0 & 6 g
}
\right)
%\eqno{(7)}
$$
Now we know that
$$
L_\bbh(g)  \  = \ \left(
\matrix
{
 \Hess^\bbp_\bbh  g   -2 g I &   * &  *  \cr
*  \ &\  0 & 0\cr
* &  0 & 0
}
\right)
%\eqno{(8)}
$$
It remains only to compute the $\bbh$-symmetric part of the *'s.

Consider a vector $(a,b) \in \ch\oplus(\cv\oplus \bbr x) = \ch\oplus \ch^\perp$. We want to look
at the term $\bra {L_\bbh a} b$.  This is 
$$
\bra {L_\bbh a} b \ =\ \smfrac 14  \left\{   
\bra {L a} b + \bra {L Ia} {Ib}  + \bra {L Ja}{ Jb}  + \bra {L Ka} {Kb} .
\right\}
$$
This is the  trace of a quadratic form on $\bbh=\bbr^4$, and it can be expressed with respect
to any orthonormal basis.  As a result we may assume that $b=x \equiv e$ and $a\in\cv$.
$$
\bra {La} b = \bra {-3(Dg) a} e \ =\ -3 \bra {(Dg)a} e
$$
and 
$$
\bra {LIa} {Ie} =\bra {Ia} {LIe} = \bra {Ia} { -\mu^*_I(Dg)} \ =\  -\bra {I^2a}{(Dg)} \ =\ \bra {a}{Dg}
$$
Similarly,
$$
\bra {LJa} {Je} =\bra {LKa} {Ke} = \bra {a}{Dg}.
$$
The sum is zero. Hence we have that with respect to the decomposition $\ch\oplus \ch^\perp$
$$
|x|^2 \left ( D^2 u\right)_\bbh   \  = \ \left(
\matrix
{
\Hess^\bbp_\bbh g   -2 g I &    0  \cr
\ &\ \cr
0  &   0
}
\right).
%\eqno{(9)}
$$
The (1,1)-term is the pull-back of the quaternionic hermitian symmetric part of the 
hessian on quaternionic projective space. This completes the proof.\qed

\vfill\eject

%%%%%%%%%%%%%%%%%%%%%%%%%%%%%%%%%%%%%%%%%%%%%%
%%%%%%%%%%%%%%%%%%%%%%%%%%%%%%%%%%%%%%%%%%%%%%
%%%%%%%%%%%%%%%%%%%%%%%%%%%%%%%%%%%%%%%%%%%%%%
%%%%%%%%%%%%%%%%%%%%%%%%%%%%%%%%%%%%%%%%%%%%%%
%%%%%%%%%%%%%%%%%%%%%%%%%%%%%%%%%%%%%%%%%%%%%%

\centerline{\headfont  Appendix A.  Further Discussion of  Examples.}
 \medskip 

In this appendix we examine specific subequations of Riesz characteristic $p$,
$1\leq p<\infty$ in more detail. We consider two types:  cone subequations 
and convex cone subequations, and  in both cases the subequations will always be
\ST.  It may be of some surprise that in each of these two  categories there is 
a unique largest and smallest subequation.

\bigskip

\centerline
{
\bf The Largest / Smallest Characteristic $p$ Subequation.
}
\medskip
We first consider the category of cone  subequations.
 For $A\in\Symn$ let 
$\l_1(A) \leq\cdots\leq \l_n(A)$ denote the ordered eigenvalues of $A$, and 
set $\l_{\rm min}(A) \equiv \l_1(A)$ and $\l_{\rm max}(A)\equiv \l_n(A)$.  We then
define
$$
F^{\rm min/max}_p \ \equiv \  \left\{A : \l_{\rm min}(A)  + (p-1) \l_{\rm max}(A) \ \geq\ 0 \right \}
\eqno{(A.1)}
$$
$$
F^{\rm min/2}_p \ \equiv \  \left\{A : \l_{\rm min}(A)  + (p-1) \l_{2}(A) \ \geq\ 0 \right \}
\eqno{(A.2)}
$$
It is clear from Definition \CC.4 that both of these subequations has Riesz characteristic $p$. 
These are the largest and smallest cone subequations with this property.

\Lemma {A.1}  
{\sl
Let $F$ be an \ST cone subequation of Riesz characteristic $p$. Then}
$$
F^{\rm min/2}_p\ \ss\ F\ \ss\ F^{\rm min/max}_p.
$$

\medskip
\noindent
{\bf Note.}  One computes that the dual of this largest subequation $F^{\rm min/max}_p$ is $F^{\rm min/max}_q$ where
$(p-1)(q-1)=1$.  Compare this with (\CC.15) in Part I which says that $(p_F-1)(q_F-1)\geq 1$
for any subequation $F$.  Also see Example \JJJ.14 in Part I.

\medskip

\noindent
{\bf Proof of Lemma A.1.}  Each $A\in \Symn$ can be written as a sum $A = \l_1 P_{e_1} + \cdots + \l_n P_{e_n}$
using the ordered eigenvalues of $A$. Set $B_0 \equiv \l_1 P_{e_1} + \l_2 P_{e_1^\perp}$,
and $B_1 \equiv \l_1 P_{e_1} + \l_n P_{e_1^\perp}$, 
and note that $B_0\leq A \leq B_1$.

If $A\in F^{\rm min/2}_p$, then $\l_1 + (p-1)\l_2\geq0$. Thus, $B_0 \in F^{\rm min/2}_p$.
Since $F^{\rm min/2}_p$ and $F$ have the same increasing radial profile $E^\uparrow$
given by (\CC.1) in Part I  (and $\l_2\geq0$), we conclude that $B_0\in F$.  
However, $B_0\leq A$ proving that $A\in F$.

For the other inclusion, pick $A \in F$.  Since $F\ss \cpt$, we have $\l_{\rm max} \geq0$.
Now $A\leq B_1$ implies $B_1\in F$. Again $F$ and $F^{\rm min/max}_p$ 
have the same increasing radial profile $E^\uparrow$ given by (\CC.1) in Part I.  Therefore,
$B_1 \in F^{\rm min/max}_p$.  This implies by definition that $A \in F^{\rm min/max}_p$.\qed

 The largest and smallest characteristic $p$ subequations  in the convex cone case are different in dimensions $\geq3$ (see Section \HH\  in Part I for the definitions of $\cp_p$ and $\cp(\d)$).

\Lemma {A.2}
{\sl
Let $F$ be an O$(n)$-invariant convex cone subequation of Riesz characteristic $p$. Then $1\leq p\leq n$ and 
$$
\qquad\qquad\qquad \cp_p\ \ss\ F\ \ss\ \cp(\d) 
 \qquad\quad {\rm where \ \ \ }  \d\ =\ {(p-1)n \over n-p}.
$$
and so the Riesz characteristic of $\cp(\d)$ is $p$.}
\pf
The first inclusion follows from the fact that $-(p-1) P_e + P_{e^\perp}$ 
generate the extreme rays in $\cp_p$.
This is proved in  [\HLPCON, Thm. 5.1c].
The second inclusion  is Proposition \JJJ.9 in Part I.\qed

\vskip .3in

\centerline{\bf O(n)-Invariant Subequations.}\medskip

Such a subequation $F$ determines a subset $E\ss\rn$ consisting of the $n$-tuples
$(\l_1(A),...,\l_n(A))$ of eigenvalues of $A$.  Consider $\l(A) = (\l_1(A),...,\l_n(A))$
as a multi-valued map $\l:\Symn \to \rn$.  Then we define $E\equiv  \l(F)$.  The set
$E$ is closed and symmetric (invariant under the permutation of coordinates in $\rn$).
In addition,
$$
E \ \ {\rm is}\ \ \rn_+\ {\rm (positive \ orphant)  \  monotone, \ i.e., } \ \ E+  \rn_+ = E, 
\eqno{(A.3)}
$$
since the ordered eigenvalues are $\cp$-monotone.

\Def{A.3}  A closed symmetric subset $E\ss\rn$ (with $\emptyset\neq E \neq \rn$)
will be called a {\bf universal eigenvalue subequation} if $E$ is $\rn_+$-monotone.
\medskip

Note that this is an abuse of language since $E$ itself is not a subequation.
The ``universal'' nature of $E$ will be described later.  However, such a set $E$ determines the 
O$(n)$-invariant subequation 
$$
F\ =\ \l^{-1}(E).
\eqno{(A.4)}
$$
Note that
$$
F\ \ {\rm is\ a\ cone\ }\qquad\iff\qquad E\ \ {\rm is\ a\ cone, \  and}
\eqno{(A.5)}
$$
$$
F\ \ {\rm is\ convex\ }\qquad\iff\qquad E\ \ {\rm is\ convex.}
\eqno{(A.6)}
$$
Of course, $\cp$ and $\rn_+$ correspond, i.e.,  $\cp =   \l^{-1}(\rn_+)$.
The Riesz characteristic of $F$ is easily computed from its eigenvalue profile  $E$.
\medskip\centerline
{
The increasing Riesz characteristic of $F$ equals $\sup\{p : (-(p-1), 1, ... , 1) \in E\}$.
}
\medskip\centerline
{
The decreasing Riesz characteristic of $F$ equals $\sup\{q : (-1, ..., -1, (q-1)) \in E\}$.
}
\medskip

It is also worth noting that if $E$ and $F$ correspond, the $\wt E$ and $\wt F$ correspond.

\vskip .3in

\centerline{\bf The Complex and Quaternionic Analogues of an O(n)-Invariant Subequation.}\medskip

As described in Example  \HH.7  each O$(n)$-invariant subequation $F$
on $\rn$     canonically determines a  U$(n)$-invariant subequation $F^\bbc$ on $\bbc^n$ 
and an Sp(n)-invariant subequation $F^\bbh$  on $\bbh^n$.  In both cases the subequation
is given by requiring that the  $n$ eigenvalues of the hermitian symmetric
part of the matrix lie in $E_F$. That is, if $E$ is defined by: 
$A\in F \iff \l(A) \in E$, then 
$$
A\in F^\bbc \ \ \iff\  \ \l^\bbc_k(A) \in E
\quad\and\quad
A\in F^\bbh \ \ \iff\ \  \l^\bbh_k(A) \in E
\eqno{(A.7)}
$$
where $\l^\bbc_k(A) = \l_k(A_\bbc)$ and $\l^\bbh_k(A) = \l_k(A_\bbh)$

The associated Riesz characteristics are given in Lemma \HH.9.

This classifies all the U$(n)$-invariant subequations $F$ with the property that:
$$
F \ =\ \pi_\bbc^{-1}\left(F^\bbc\right) \qquad {\rm where} \ \ \pi_\bbc(A) \ =\ A_\bbc
\eqno{(A.8)}
$$
and $F^\bbc$ is a subset of the hermitian symmetric matrices.
Similarly, it defines all the Sp$(n)$-invariant  subequations $F$ with the property that:
$$
F \ =\ \pi_\bbh^{-1}\left(F^\bbh\right) \qquad {\rm where} \ \ \pi_\bbh(A) \ =\ A_\bbh
\eqno{(A.9)}
$$
and $F^\bbh$ is a subset of the quaternionic hermitian symmetric matrices.

The largest / smallest results (Lemmas A.1 and A.2 in the $\bbr$ case)
have counterparts in the $\bbc$ and $\bbh$ cases. The precise statements and their
proofs are left to the reader.

\Ex{A.4. (Lagrangian Subharmonic)}  A notable new example of a U(n)-invariant subequation
not satisfying (A.8) comes from Lagrangian geometry, namely,   the geometrically defined subequation
$F({\rm LAG})$ where $\LAG \ss G^\bbr(n, \bbc^n)$ is the set of Lagrangian $n$-planes
in $\bbc^n$.  The eigenvalues of the skew-hermitian part $A_\bbc^{\rm skew}$ come in pairs
$\l_1, -  \l_1, \l_2, -\l_2, ... , \l_n, -\l_n$.  The subequation $F({\rm LAG})$ is determined by
a constraint on these eigenvalues together with the real  trace  $t=\tr(A)$, namely
$$
{t\over 2} \ \pm \ \l_1 \ \pm \ \l_2\ \pm\ \cdots\ \pm\ \l_n\ \geq\ 0
\eqno{(A.10)}
$$
for all $2^n$ choices of $\pm$.  There is a polynomial operator $M$ on $\Sym_\bbr(\bbc^n)$
analogous to the determinants $\det_\bbr A$, $\det_\bbc A$ and $\det_\bbh A$ 
(defined as the product of the eigenvalues in all cases), namely
$$
M_{\LAG}(A) \ =\ \prod_{2^n\ { \rm times}}  \left({t\over 2} \ \pm \ \l_1 \ \pm \ \l_2\ \pm\ \cdots\ \pm\ \l_n\right )
\eqno{(A.11)}
$$
(see [\HLDD, page 433]). Of course since $\LAG\ss G(n,\bbr^{2n})$ and $F$ is geometrically defined 
by $\LAG$, $F$ has Riesz characteristic $n$.

\def\ISO{{\rm ISO}}

\Ex{A.4$_p$. (Isotropic Subharmonic)} The previous example can be generalized as follows.
For each integer $p$, $1\leq p\leq n$ we consider the set 
$$
\ISO_p \ =\ \{W\in G^\bbr(p, \bbc^n) : W \ {\rm is\  an\  isotropic\ } p\ {\rm plane}\}.
$$
Recall that a real $p$-plane $W$ in $\bbc^n=\bbr^{2n}$ is {\sl isotropic} if 
$$
v\perp Jw\qquad \forall\, v,w\in W,
$$
i.e., the K\"ahler form $\o$ satisfies $\o\bigr|_W=0$. Note that $\ISO_n = \LAG$ and $\ISO_1=\cp$.
For all $p>1$, the set  $\ISO_p$ has the transitivity property, and so Theorem \KK.2 applies.

 Associated to this set is the subequation
$$
\eqalign
{
F(\ISO_p)\ 
                   &=\   \{ A\in\Sym(\bbr^{2n}) : \forall\, W\in \ISO_p,\  \tr_W A\geq 0\}    \cr
                 &=\   \{ A\in\Sym(\bbr^{2n}) : \smfrac {p}{n} t-\l_{n-p+1}-\cdots-\l_n\geq 0\}    \cr 
}
$$
where $0\leq \l_1\leq \cdots \leq \l_n$ are  as in Example A.4. The dual subequation is
$$
\eqalign
{
\wt F(\ISO_p)\ 
                   &=\   \{ A\in\Sym(\bbr^{2n}) : \exists \, W\in \ISO_p,\  \tr_W A\geq 0\}    \cr
                 &=\   \{ A\in\Sym(\bbr^{2n}) : \smfrac {p}{n} t+\l_{n-p+1}+\cdots+\l_n\geq 0\}    \cr 
}
$$
Associated to this problem we have the polynomial
$$
{\bf M}_{\ISO_p}(A) \ =\  \prod_{|I|=p\,{\rm and\ }\pm}
 \left(\smfrac {p}{n} t \pm\l_{i_1} \pm\cdots\pm \l_{i_p}\right)
$$
which is   a factor of $\det(D_{A_{LAG}})$ where 
$D_{A_{LAG}} : \L^p \bbr^{2n}\to \L^p \bbr^{2n}$ 
is the extension of $A_{LAG}$ as a derivation.
As above we have that any $C^2$ function $u$ which is 
$\ISO_p$-harmonic satisfies the differential equation
$$
{\bf M}_{\ISO_p}(\Hess\, u) \ =\ 0
$$

\vskip .3in

\centerline{\bf Subequations Arising from G\aa rding Operators.}\medskip

G\aa rding's beautiful theory of hyperbolic polynomials provides a surprisingly rich collection
of nonlinear operators. (This connection is mentioned in Krylov [Kr]).  Moreover, associated 
with each such `` G\aa rding operator'' there are many actual subequations.  Here we provide
 a brief overview. We first start  with an operator and   discuss how the many associated subequations
 are universally constructed. We then describe  three of the basic ways of constructing new 
G\aa rding operators from a given one.  These repeatable processes lead to a vast array of 
G\aa rding operators, starting with just one. 

See for example [\HLHPP]  for a self-contained development of G\aa rding's theory.
His two fundamental results can be summarized by saying that:
$$
{\rm The \ Garding\ cone\ } \Gamma\ {\rm is\ convex,\ and}
\eqno{(A.12)}
$$
$$
{\rm The \ Garding\ eigenvalues\  are\ } \Gamma-{\rm monotone.}
\eqno{(A.13)}
$$

\Def{A.5}  A homogeneous real polynomial $M$ of degree $m$ on the space $\Symn$ of 
second derivatives, with $M(I)>0$, is a {\bf G\aa rding operator} if:
\medskip

(1) For each $A\in\Symn$ the polynomial $M(sI+A)$ has $m$ real roots ($M$ is $I$-{\sl hyperbolic}), and

\medskip

(2)  The G\aa rding cone $\G$, defined as the connected component of $I$ in $\{M(A)>0\}$, satisfies 
{\sl positivity} $\G+\cp\ss\G$.
\medskip

The primary subequation associated with the G\aa rding operator $M$ is the closure of the
G\aa rding cone $\G$, which is a convex cone subequation. However, there are many others.

\def\l{\lambda^M}

\Def{A.6} The negatives of the roots  of $M(sI+A)=0$ are called the $M$-{\sl eigenvalues of $A$},
and are denoted by $\l(A) \equiv (\l_1(A), ... , \l_m(A))$.  They are well defined up to permutations.
Note that $M(A) = \l_1(A) \cdots \l_n(A)$.

\Def{A.7} The {\bf ${\bf  k}^{\rm th}$ branch of the equation $M(A)=0$} is defined to be the set
$$
\{\l_k(A)\ \geq\ 0\}
$$
where $\l_1(A) \leq \cdots  \leq \l_m(A)$ are the {\sl ordered} eigenvalues of $A$.

\medskip
An important part of the theory shows that the ordered eigenvalues  are strictly $\G$-monotone.
Since $\cp\ss\overline\G$ by (A.13), we have that each of the $m$ branches of $\{M(A)=0\}$
is a subequation.  

Note that $\overline\G = \{\l_{\rm min}(A)\geq0\}$ is the primary branch, and its dual
subequation is $\{\l_{\rm max}(A)\geq0\}$, which is the largest branch.

The branches $\{\l_k(A)\geq0\}$, $k=1,...,m$, are the subequations most intimately associated
with the G\aa rding operator $M$ in that if a $C^2$-function $u$ is harmonic for one of these
subequations, then 
$$
M\left( D_x^2 u\right) \ =\ 0.
\eqno{(A.14)}
$$
However, there are many others, all constructed exactly as in the O(n)-invariant case.

Now we make full use of the concept  (Def. A.3) of a universal eigenvalue subequation.

\Prop{A.8}  {\sl
Given a  universal eigenvalue subequation $E\ss \bbr^m$, each G\aa rding operator of degree
$m$ on $\Symn$ determines a subequation on $\rn$, namely
$$
F_E \ \equiv \ \{ A : \l(A) \in E\}.
$$
This subequation is $\overline \G$-monotone (not just $\cp$-monotone).  Moreover, 
$F_E$ is a cone if and only if $E$ is a cone., and $F_E$ is convex if and only if $E$ is convex.
}
\pf This is straightforward except for the last assertion which is due to [BGLS].
\medskip

For example, $E = \bbr^m_+$ is the universal ``Monge-Amp\`ere subequation''
inducing the subequation $\overline \G$ for each degree $m$ G\aa rding operator
 $M(A)= \l_1(A) \cdots \l_m(A)$.

We complete this discussion of G\aa rding operators by describing three of the basic 
methods of constructing new  G\aa rding operators from a given  G\aa rding operator 
$M$ of degree $m$. To be specific the reader may want to start with one of the 
basic operators $\det (A_K)$ for $K=\bbr,\bbc$ or $\bbh$.

\medskip
\noindent
{\bf I. The Derived or Elementary Symmetric Operator.}  With $k=1,...,m$  we set
$$
\s_k(A) \ \equiv \ {1\over (n-k)!} {d^{m-k}\over dt^{m-k}} M(A+tI)\bigr|_{t=0}
\ =\ \sum_{i_1<\cdots i_k} \l_{i_1}(A) \cdots  \l_{i_k}(A) 
$$

\medskip
\noindent
{\bf II. The $p$-Convexity Operator.}  For each real number $p$ with $1\leq p\leq m$  set
$$
\Sigma_p(A) \ \equiv\  \prod\left\{ \l_{i_1}(A) + \cdots + \l_{i_{[p]}} (A)+ (p-[p]) \l_j(A)\right \}
$$
where the product is taken over  all  increasing multi-indicies $I=(i_1,..., i_{[p]})$ 
and all $j\notin I$.

\medskip
\noindent
{\bf III. The $\d$-Uniformly Elliptic Regularization Operator.}  With $0\leq \d \leq \infty$ (and 
renormalizing at $\d=\infty$) set
$$
M^\d(A) \ \equiv\ \prod_{j=1}^m \left\{ \l_j(A) + {\d\over n}(\tr A)\right\}
$$

\def\l{\lambda}

\Remark{A.9} The first process lowers the degree of the operator.  The second process
raises the degree of the operator, and the degree remains the same in the third process.
Each of these construction can be applied repeatedly. In the third process  nothing
new is obtained from this iteration.  However,  iterating the first  process  produces a finite
number of new G\aa rding operators, and iterating the second  produces an 
infinite family of new ones.  Moreover, one can apply any sequence of the three 
operations, thereby producing a huge collection of G\aa rding operators all dependent
on the primary operator.

\vskip .3in

\centerline{\bf Elliptic Regularization  -- Subequation Expansion/Contraction.}\medskip

For each $r>0$ consider the linear map
$$
\Phi_r(A)\ \equiv\ rA + (1-r)(\tr\,A)\smfrac 1n I 
\ =\  r\left (A - (\tr \,A)\smfrac 1n  I\right) +  (\tr \,A)\smfrac 1n  I.
\eqno{(A.15)}
$$
The restriction of $\Phi_r$ to each affine hyperplane $\{\tr\,A=\l\}$ is the $r$-homothety
(multiplication by $r$) about the center ${\l\over n}I$.  This follows from the second equality.
The inverse is
$$
\Phi_{1\over r} \ =\ {1\over r} \left ( A + (r-1) (\tr\,A)\smfrac 1n I \right).
\eqno{(A.16)}
$$

\Def{A.10} Suppose $\d \equiv r-1\geq0$ and $F$ is a cone subequation.
Then 
$$
F(\d) \ =\ \Phi_r(F) \ =\ \{A : r\Phi_{1\over r} (A) = A+ \d(\tr\,A)\smfrac 1n I\in F\}
$$
is called the $\bf r^{\bf th}$ {\bf expansion of $F$}.
\medskip

Note that $F(\d)$ is also a subequation for all $\d>0$ since the homothety factor $r\geq1$.
Note also that if $F$ is a convex cone contained in $\Int \D = \{\tr\,A> 0\}$,
then $F(\d)$ ranges from $F$ to $\D$ as $\d$ ranges from 0 to $\infty$.
Finally, note that $\partial F(\d) = \Phi_r(\partial F)$.

\Prop{A.11} {\sl
Suppose $F$ is a cone subequation with (Riesz) characteristic $p=p_F$ and dual
 (Riesz) characteristic $q=q_F$.  Then the $\d$-uniformly elliptic cone subequation 
 $F(\d) = \Phi_r(F)$ ($\d\equiv r-1\geq0$)
 has its two characteristics given by the same function
 $$
 \eqalign
 {
p_{F(\d)}\ &= \ {n(1+\d)p  \over n+\d p} \ =\ p + {\d p (n-p) \over n + \d p}   \cr
q_{F(\d)}\ &= \ {n(1+\d)q  \over n+\d q} \ =\ q + {\d q (n-q) \over n + \d q}   \cr
 }
\eqno{(A.17)}
$$
 These formulas hold when $p_F=\infty$ or when  $q_F=\infty$, that is}
 $$
  \eqalign
 {
p_{F}\ &=\ \infty\quad\Rightarrow\quad p_{F(\d)} = {n(1+\d) \over \d} \cr
q_{F}\ &=\ \infty\quad\Rightarrow\quad  q_{F(\d)} = {n(1+\d) \over \d} \cr
}
 \eqno{(A.18)}
$$

\pf
Note that $A\equiv P_{e^\perp} -(p-1) P_e \in\partial F \iff \Phi_r(A) \in \partial \Phi_r(F) = \partial F(\d)$ and
$$
\Phi_r(A) \ =\ (1+\d) P_{e^\perp} - (1+\d) (p-1) P_e - {\d(n-p)\over n} I
\ =\ {n+\d p \over n} \left[ P_{e^\perp}  - \left(  {n(1+\d)p  \over n+\d p} -1   \right) P_e    \right].
$$
Finally, since $-A\in \partial F \iff -\Phi_r(A) \in \partial \Phi_r(F)= \partial F(\d)$,
the formula for $q_{F(\d)}$ as a function of $q_F$ is the same as the
formula for  $p_{F(\d)}$ as a function of $p_F$. \qed

\Prop{A.12}  {\sl
If $F$ is $M$-monotone, then $F(\d)$ is $M(\d)$-monotone.}

\pf Straightforward.

\Ex{A.13}  As $\d$ ranges from 0 to $\infty$, $\cp(\d)$ increases from $\cp$ to $\D$.
Each $\cp(\d)$ is a convex cone ; and with $\d>0$ small, these subequations form a 
``fundamental system'' of conical neighborhoods of $\cp$.  Consequently,
they provide one of the nicer definitions of uniform ellipticity.  Namely, a subequation
$F$ is {\bf $\d$-uniformly elliptic} if
$$
F+\cp(\d) \ \ss\ F.
\eqno{(A.19)}
$$

Since $F(\d) + \cp(\d) \ss F(\d)$, each $F(\d)$ is automatically $\d$-uniformly elliptic.
For this reason, $F(\d)$ is also called the {\bf $\d$-elliptic regularization of $F$}
(cf. [Kr]).

\vfill \eject

\centerline{\bf REFERENCES}

\vskip .2in

\noindent
\item{[A$_1$]}   S. Alesker,  {\sl  Non-commutative linear algebra and  plurisubharmonic functions  of quaternionic variables}, Bull.  Sci.  Math., {\bf 127} (2003), 1-35. also ArXiv:math.CV/0104209.  

\smallskip

\noindent
\item{[A$_2$]}   \ \----------,   {\sl  Quaternionic Monge-Amp\`ere equations}, 
J. Geom. Anal., {\bf 13} (2003),  205-238.

 ArXiv:math.CV/0208805.  

\smallskip

\noindent
\item{[AV]}    S. Alesker and M. Verbitsky,  {\sl  Plurisubharmonic functions  on hypercomplex manifolds and HKT-geometry},  J. Geom. Anal. {\bf 16} (2006), no. 3, 375Ð399.

\smallskip

\noindent
\item{[BGLS]}   H. Bauschke, O. G\"uler,  A. Lewis and H. Sendov, {\sl   Hyperbolic polynomials and convex analysis},
Canad. J. Math. {\bf 53} no. 3, (2001), 470-488. 

 \smallskip

 \noindent 
\item {[\HLCG]}   F. R. Harvey and H. B. Lawson, Jr,  {\sl Calibrated geometries}, Acta Mathematica 
{\bf 148} (1982), 47-157.

 \smallskip

\item {[\HLPCG]}  \ \----------, 
 {\sl  An introduction to potential theory in calibrated geometry}, Amer. J. Math.  {\bf 131} no. 4 (2009), 893-944.  ArXiv:math.0710.3920.

\smallskip

\item {[\HLDD]}   \ \----------,  {\sl  Dirichlet duality and the non-linear Dirichlet problem},    Comm. on Pure and Applied Math. {\bf 62} (2009), 396-443. ArXiv:math.0710.3991

\smallskip

\item {[\HLDDR]}  \ \----------,   {\sl  Dirichlet duality and the nonlinear Dirichlet problem
on Riemannian manifolds},  J. Diff. Geom. {\bf 88} (2011), 395-482.   ArXiv:0912.5220.
\smallskip

\item {[\HLHPP]}  \ \----------,   {\sl  G\aa rding's theory of hyperbolic polynomials},  
 Comm. Pure Appl. Math. {\bf 66} (2013), no. 7, 1102Ð1128.

\smallskip

\item {[\HLREST]}  \ \----------, {\sl  The restriction theorem for fully nonlinear subequations}, 
   Ann. Inst.  Fourier (to appear). 
ArXiv:1101.4850.

\smallskip

\item {[\HLPCON]}   \ \----------,  {\sl  p-convexity, p-plurisubharmonicity  and the Levi problem },
   Indiana Univ. Math. J.  {\bf 62} No. 1 (2014), 149-170.  ArXiv:1111.3895.

\smallskip

\item {[\HLSURVEY]}   \ \----------,  {\sl  Existence, uniqueness and removable singularities
for nonlinear partial differential equations in geometry},  pp. 102-156 in ``Surveys in Differential Geometry 2013'', vol. 18,  
H.-D. Cao and S.-T. Yau eds., International Press, Somerville, MA, 2013.
ArXiv:1303.1117.

\smallskip

\item  {[\HLRS]} \ \----------, {\sl  Removable singularities for nonlinear subequations}, 
 ArXiv:1303.0437.
\smallskip

\noindent
\item  {[\HLTANI]} \ \----------, {\sl Tangents to subsolutions  --  existence and uniqueness, I}, 
(to appear).

\smallskip

 \noindent
\item{[K]}
C. Kiselman, {\sl Tangents of plurisubharmonic functions}, 
 International Symposium in Memory of Hua Loo Keng, Vol. II (Beijing, 1988), 157Ð167, Springer, Berlin, 1991. 

\smallskip

   \noindent
\item{[Kr]}    N. V. Krylov,    {\sl  On the general notion of fully nonlinear second-order elliptic equations},    Trans. Amer. Math. Soc. (3)
 {\bf  347}  (1979), 30-34.

\smallskip

 \item {[R]} 
R, T.   Rockafellar, Convex analysis. Princeton Mathematical Series, No. 28,  Princeton University Press, Princeton, N.J.,  1970.

\vfill\eject

\end

\vskip .3in
 
\centerline{\bf ABSTRACT} \medskip
 % \font\abstractfont=cmr10 at 10 pt

  {{\parindent= .44in
\narrower  \noindent
This part II of the paper is concerned with questions of existence
and uniqueness of tangents in the special case of  $\GG$-plurisubharmonic functions,
where $\GG \ss G(p,\rn)$ is a  compact subset of the 
Grassmannian of $p$-planes in $\rn$, invariant under a subgroup $K\ss {\rm O}(n)$
which acts transitively on $S^{n-1}$.  An u.s.c. function $u$ on an
open set $\O\ss\rn$ is $\GG$-plurisubharmonic if its restriction to
$\O\cap W$ is subharmonic for every affine $\GG$-plane $W$.
Tangents to $u$ at a point $\point$ are the cluster points of $u$ under 
a natural flow (or blow-up) at $\point$. 
Tangents always exist and are $\GG$-harmonic at all points of continuity.
A homogeneity property is proved for all tangents in these geometric cases.
A major result concerns the Strong Uniqueness of Tangents -- a form of regularity
which has  consequences  for the density function established  in Part I.
It means that all tangents are unique and of the form $\den K_p$ where
$K_p$ is the Riesz kernel and $\den$ is the density of $u$ at the point.
Strong uniqueness is shown to hold for all but
a small handful of cases when $K= {\rm O}(n), {\rm U]}(n)$ or Sp$(n)$.
It also holds for essentially  all interesting $\GG$ which arise in calibrated geometry.

In the cases corresponding to the real, complex and quaternionic Monge-Amp\`ere
equations (convex functions, and complex and quaternionic plurisubharmonic
functions) tangents, which are far from unique, are systematically studied and
classified.

}}

\vskip .3in
 
\centerline{\bf ABSTRACT} \medskip
 % \font\abstractfont=cmr10 at 10 pt

  {{\parindent= .44in
\narrower  \noindent
This part II of the paper is concerned with questions of existence
and uniqueness of tangents in the special case of  $\GG$-plurisubharmonic functions,
where $\GG \ss G(p,\rn)$ is any compact $K$-invariant subset of the 
Grassmannian of $p$-planes in $\rn$, where $K\ss {\rm O}(n)$
acts transitively on $S^{n-1}$.  An u.s.c. function $u$ on an
open set $\O\ss\rn$ is $\GG$-plurisubharmonic if its restriction to
$\O\cap W$ is subharmonic for every affine $\GG$-plane $W$.
A homogeneity property is proved for all tangents in these geometric cases.
Strong uniqueness of tangents is established in a wide range of cases.
In the cases corresponding to the real, complex and quaternionic Monge-Amp\`ere
equations (convex functions, and complex and quaternionic plurisubharmonic
functions) tangents, which are far from unique, are systematically studied and
classified.

}}

\vskip .3in
 
\centerline{\bf ABSTRACT} \medskip
 % \font\abstractfont=cmr10 at 10 pt

  {{\parindent= .44in
\narrower  \noindent
This part II of the paper is concerned with questions of existence
and uniqueness of tangents in the special case of  $\GG$-plurisubharmonic functions,
where $\GG \ss G(p,\rn)$ is an arbitrary compact subset of the 
Grassmannian of $p$-planes in $\rn$.  An u.s.c. function $u$ on an
open set $\O\ss\rn$ is $\GG$-plurisubharmonic if its restriction to
$\O\cap W$ is subharmonic for every affine $\GG$-plane $W$.
A homogeneity property is proved for all tangents in these geometric cases.
Strong uniqueness of tangents is established in a wide range of cases.
In the cases corresponding to the real, complex and quaternionic Monge-Amp\`ere
equations (convex functions, and complex and quaternionic plurisubharmonic
functions) tangents, which are far from unique, are systematically studied and
classified.

}}

When $F$ is convex, there are other densities defined via the {\bf area and 
volume averages:}
$$
S(u,r) \ =\  \intave S u 
\and
V(u,r) \ =\ \intave B u
\eqno{(\AAA.6)}
$$
where $B = \{|x|\leq1\}$ and $S=\partial B$.